\documentclass{article}

\usepackage{amsmath, amscd, amsthm, amssymb, graphics}

\newcommand{\Proj}{{\rm Proj}}
\newcommand{\Spec}{{\rm Spec}}
\newcommand{\diag}{{\rm diag}}

\input amssym.def
\input amssym.tex
\newcommand{\nc}{\newcommand}
\nc{\bla}{\phantom{bbbbb}}

\newcommand{\re}{{\rm re} \, }

\newtheorem{theorem}{Theorem}[section]
\newtheorem{corollary}[theorem]{Corollary}

\newtheorem{prop}[theorem]{Proposition}
\newtheorem{definition}[theorem]{Definition}
\newtheorem{remit}[theorem]{Remark}

\newtheorem{exit}[theorem]{Example}

\newenvironment{rem}{\begin{remit}\rm}{\end{remit}}
\newenvironment{ex}{\begin{exit}\rm}{\end{exit}}
\newenvironment{defn}{\begin{definition}\rm}{\end{definition}}

\newcommand{\RR}{{\Bbb R }}
\newcommand{\CC}{{\Bbb C }}
\newcommand{\ZZ}{{\Bbb Z }}
\newcommand{\PP}{ {\Bbb P } }
\newcommand{\QQ}{{\Bbb Q }}

\newcommand{\calh}{{\mbox{$\cal H$}}}
\newcommand{\call}{{\mbox{$\cal L$}}}
\newcommand{\calo}{{\mbox{$\cal O$}}}
\newcommand{\calw}{{\mbox{$\cal W$}}}
\newcommand{\calx}{{\mbox{$\cal X$}}}
\newcommand{\caly}{{\mbox{$\cal Y$}}}

\nc{\conv}{{\rm Conv}}
\nc{\umax}{{U_{\max}}}
\newcommand{\liek}{{\bf k}}
\newcommand{\lieu}{{\bf u}}

\newcommand{\lieks}{{\bf k}^*}
\newcommand{\xg}{X/\!/G}
\newcommand{\xu}{X/\!/U}

\newcommand{\txg}{\tilde{X}/\!/G}
\newcommand{\cplusr}{(\CC^+)^r}
\newcommand{\slrplus}{SL(r+1;\CC)}
\newcommand{\glr}{GL(r;\CC)}
\newcommand{\prrplus}{{\PP^{r(r+1)}}}
\newcommand{\prr}{{\PP^{r^2}}}

\nc{\lieg}{{\bf g}}
\nc{\liegs}{{\bf g}^*}

\def\a{\alpha}
\def\b{\beta}
\def\g{\gamma}

\def\e{\epsilon}
\def\c{\chi}

\title{Quotients by non-reductive algebraic group actions}

\author{Frances Kirwan\\Mathematical Institute, Oxford OX1 3BJ, UK\\kirwan@maths.ox.ac.uk}

\begin{document}

\maketitle

This paper is dedicated to Peter Newstead, from whose Tata Institute Lecture
Notes \cite{New} I learnt about GIT and moduli spaces some
decades  ago, with much appreciation
 for all his help
and support over the years since then.

\section{Introduction}

Geometric invariant theory (GIT) was developed in the 1960s by Mumford in order to construct
quotients of  
reductive group actions
on algebraic varieties and hence to construct and study a number of moduli 
spaces, including, for example, moduli spaces of bundles over a nonsingular
projective curve \cite{GIT,New,New2}. 
Moduli spaces often arise naturally as quotients of varieties by algebraic
group actions, but the groups involved are not always reductive. For example, in the case of moduli spaces of hypersurfaces (or, more generally, complete intersections) in
toric varieties (or, more generally, spherical varieties), the group actions which arise naturally are actions of the automorphism groups of the varieties \cite{cox,coxkatz}. These automorphism groups
are not in general reductive, and when they are not reductive we cannot use classical GIT to construct
 (projective completions of) such moduli spaces as quotients for these actions.
 
 In \cite{DK} (following earlier work including \cite{F2,F1,GP1,GP2,W} and references therein)
a study was made of ways in which GIT might be generalised to non-reductive
 group actions; some more recent developments and applications can be found in
 \cite{AD,AD2}. Since every affine algebraic group $H$ has a unipotent radical $U \unlhd
 H$ such that $H/U$ is reductive, \cite{DK} concentrates on unipotent actions. It is
 shown that when a unipotent group $U$ acts linearly (with respect to an ample
 line bundle $\call$) on a complex projective variety $X$, then $X$ has invariant
 open subsets $X^s \subseteq X^{ss}$, consisting of the \lq stable' and \lq semistable'
 points for the action, such that $X^s$ has a 
 geometric quotient $X^s/U$ and $X^{ss}$ has a 
 canonical \lq enveloping quotient'  $X^{ss} \to X/\!/U$ which restricts to $X^s \to
 X^s/U$ where $X^s/U$ is an open subset of $X/\!/U$. (When it is necessary to distinguish between
 stability and semistability for different group actions on 
 $X$ we shall denote $X^s$ and $X^{ss}$ by $X^{s,U}$
 and $X^{ss,U}$.) However, in contrast to the reductive case, the natural morphism
 $X^{ss} \to X/\!/U$ is not necessarily surjective; indeed
 its image is not necessarily a subvariety of $X/\!/U$, so we do not in general
 obtain a categorical quotient of $X^{ss}$. Moreover $X/\!/U$ is in general only 
 quasi-projective, not projective, though when the ring of invariants
 $\hat{\calo}_L(X)^U = \bigoplus_{k \geq 0} H^0(X, L^{\otimes k})^U$
is finitely generated as a $\CC$-algebra then
$X/\!/U$ is the projective variety $\Proj({\hat{\calo}}_L(X)^U)$.

In order to construct a projective completion of the enveloping quotient
$X/\!/U$ when the ring of invariants ${\hat{\calo}}_L(X)^U$ is not
finitely generated, and to understand its geometry, it is convenient 
to transfer the problem of constructing a quotient for the
$U$-action to the construction of a quotient for an action of a reductive
group $G$ which contains $U$ as a subgroup, by finding a \lq reductive
envelope'. This is a projective completion 
$$\overline{G \times_U X}$$
of the quasi-projective variety $G \times_U X$ (that is, the quotient of $G \times X$ by the free action of $U$ acting diagonally on the left on $X$
and by right multiplication on $G$), with a linear $G$-action on 
$\overline{G \times_U X}$ which restricts to the induced $G$-action
on $G \times_U X$, such that \lq sufficiently many' $U$-invariants on
$X$ extend to $G$-invariants on $\overline{G \times_U X}$.
If (as is always possible) we choose the linearisation on
$\overline{G \times_U X}$ to be ample (or more generally to
be \lq fine'), then the classical GIT quotient
$$\overline{G \times_U X}/\!/G$$
is a (not necessarily canonical) projective completion
$\overline{X/\!/U}$ of $X/\!/U$, and hence also of its
open subset $X^s/U$ if $X^s \neq \emptyset$, and we
have
$$X^{\bar{s}} \subseteq X^s \subseteq X^{ss} \subseteq X^{\bar{ss}}$$
where $X^{\bar{s}}$ (respectively $X^{\bar{ss}}$) denotes the
open subset of $X$ consisting of points of $X$ which are stable
(respectively semistable) for the $G$-action on
$\overline{G \times_U X}$ under the inclusion 
$$X \hookrightarrow G \times_U X \hookrightarrow \overline{G \times_U X}.$$
In principle, at least, we can apply the methods of classical
GIT to study the geometry of this projective completion
$\overline{X/\!/U}$ in terms of the geometry of the $G$-action on 
$\overline{G \times_U X}$; for example, techniques
from symplectic geometry can be used to study the topology
of GIT quotients of 
complex projective varieties by  complex reductive group actions \cite{AB,JK,JKKW,K,K2,K3,K4}.

Given a linear $H$-action on $X$ where $H$ is an affine algebraic
group with unipotent radical $U$, if an ample reductive envelope
$\overline{G \times_U X}$ is chosen in a sufficiently
canonical way that the GIT quotient 
$\overline{G \times_U X}/\!/G = \overline{X/\!/U}$
inherits an induced linear action of the reductive group
$R= H/U$, then 
$$(\overline{X/\!/U})/\!/R$$
is a projective completion of the geometric quotient
$$X^{s,H}/H = (X^{s,H}/U)/R$$
where $X^{s,H}$ is the inverse image in $X^{s,U}$ of the
open set of points in 
$X^{s,U}/U \subseteq X/\!/U \subseteq \overline{X/\!/U}$
which are stable for the action of $R=H/U$. Moreover we can
study the geometry of $(\overline{X/\!/U})/\!/R$ in terms
of that of $\overline{X/\!/U}$ using classical GIT and symplectic geometry.

The aim of this paper is to discuss, for suitable actions
on projective varieties $X$
of a non-reductive affine algebraic group $H$ with unipotent radical $U$,
how to choose a reductive group $G \geq U$ and reductive
envelopes $\overline{G \times_U X}$. 
In particular we will study the family of examples given by moduli spaces of hypersurfaces
in the weighted projective plane $\PP(1,1,2)$, obtained as quotients
by linear actions of the automorphism group $H$ of $\PP(1,1,2)$ (cf. \cite{F3}). This
automorphism group
is a semidirect product $H= R \ltimes U$ where $R \cong GL(2;\CC)$ is reductive
and $U \cong (\CC^+)^3$ is the unipotent radical of $H$, acting
on $\PP(1,1,2)$ as
$$[x:y:z] \mapsto [x:y:z+ \lambda x^2 + \mu xy + \nu y^2]$$
for $\lambda, \mu, \nu \in \CC$.

For simplicity we will work over $\CC$ throughout. The layout of the paper is as follows. 
$\S$2 gives a very brief review of classical GIT 
for reductive group actions, while $\S$3 describes the results of
\cite{DK} on non-reductive actions and the construction of reductive
envelopes. $\S$4 discusses the choice of reductive envelopes for actions
of unipotent groups of the form $\cplusr$. Finally $\S$5 considers the case of the automorphism
group of $\PP(1,1,2)$ acting on spaces of hypersurfaces. 

This paper is based heavily on joint work with
Brent Doran \cite{DK,DK2}. I would like to thank him for many
helpful discussions, and also Keith Hannabuss for pointing me to
some very useful references.

\section{Mumford's geometric invariant theory}

In the preface to the first edition of \cite{GIT}, Mumford states
that his goal is ``to construct moduli schemes for various types of
algebraic objects'' and that this problem ``appears to be, in
essence, a special and highly non-trivial case'' of the problem of
constructing orbit spaces for algebraic group actions. More
precisely, when a family ${\cal X}$ of objects with parameter space $S$ has the
local universal property (that is, any other family is locally equivalent 
to the pullback of ${\cal X}$ along
a morphism to $S$) for a given moduli problem, and a group
acts on $S$ such that objects parametrised by points in $S$ are
equivalent if and only if the points lie in the same orbit, then the
construction of a coarse moduli space is equivalent to the
construction of an orbit space for
the action (cf. \cite[Proposition 2.13]{New}).
 Here, as in \cite{New}, by an orbit space 
we mean a $G$-invariant morphism
$\phi:S \to M$ such that every other $G$-invariant morphism
$\psi:S \to M$ factors uniquely through $\phi$ and $\phi^{-1}(m)$ is
a single $G$-orbit for each $m \in M$. 

Of course such orbit spaces do { not} in general exist,
in particular because of the jump phenomenon: there 
may be orbits contained in the
closures of other orbits, which means that the set of all orbits 
 cannot be
endowed naturally with the structure of a variety.
This is the situation with which Mumford's geometric 
invariant theory \cite{GIT} attempts to deal, when
the group acting is reductive, telling
us (in suitable circumstances) both how to throw out certain (unstable) orbits in
order to be able to construct an orbit space, and how to construct a projective
completion of this orbit space. Mumford's GIT is reviewed briefly next; for
more details see \cite{New2} in these proceedings, or \cite{Dolg,GIT,New,PopVin}.

\begin{ex} Let $G=SL(2;\CC)$ act on $(\PP^1)^4$ in the standard way.
Then
$$\{ (x_1,x_2,x_3,x_4)\in (\PP^1)^4: x_1=x_2=x_3=x_4 \}$$
is a single orbit which is contained in the closure of every
other orbit. On the other hand, the open subset
of $(\PP^1)^4$ where $x_1,x_2,x_3,x_4$ are  distinct
has an orbit space which, using the cross-ratio, can be identified 
with
$\PP^1 - \{0,1,\infty\}.$
\end{ex}

\subsection{Classical geometric invariant theory}

Let $X$ be a complex projective variety and let $G$ be a complex
reductive group acting on $X$. 
 To apply geometric invariant theory
we require a  {linearisation} of the action; that is, a
line bundle $L$ on $X$ and a lift of the action of $G$ to $L$.
Usually $L$ is assumed to be ample, and then
we lose little generality in supposing that for some
projective embedding 
$X \subseteq \PP^n$ 
the action
of $G$ on $X$ extends to an action on $\PP^n$ given by a
representation 
$$\rho:G\rightarrow GL(n+1),$$
and taking for $L$ the hyperplane line bundle on $\PP^n$.
We have an induced action of $G$ on the
homogeneous coordinate ring 
$${\hat{\calo}}_L(X) = \bigoplus_{k \geq 0} H^0(X, L^{\otimes k}) $$
of $X$.
The
subring ${\hat{\calo}}_L(X)^G$ consisting of the elements of ${\hat{\calo}}_L(X)$
left invariant by $G$ is a finitely generated graded complex algebra
because $G$ is reductive \cite{GIT}, and
so we can define the GIT quotient $X/\!/G$ to be the variety  $\Proj({\hat{\calo}}_L(X)^G)$. The inclusion of
${\hat{\calo}}_L(X)^G$ in ${\hat{\calo}}_L(X)$ defines a { rational} map $q$ from $X$
to $X/\!/G$, but because there may be points of $X \subseteq \PP^n$ where
every $G$-invariant polynomial vanishes this map will not in general be
well-defined everywhere on $X$. 
The set $X^{ss}$ of {\em semistable} points in $X$ is
the set of those $x \in X$ for which there exists
some $f \in {\hat{\calo}}_L(X)^G$ not vanishing at $x$. Then the rational
map $q$ restricts to a surjective $G$-invariant morphism from
the open subset $X^{ss}$ of $X$ to the quotient variety $X/\!/G$. However
$q:X^{ss} \to X/\!/G$ is still not in general an orbit
space: when $x$ and $y$ are semistable points of $X$ we have
$q(x)=q(y)$ if and only if the closures $\overline{O_G(x)}$ and
$\overline{O_G(y)}$ of the $G$-orbits of $x$ and $y$ meet in
$X^{ss}$. 

A {\em stable} point of $X$ (\lq properly stable' in
the terminology of \cite{GIT})
is a point $x$
of $X^{ss}$ with a neighbourhood in $X^{ss}$ such that every $G$-orbit 
meeting this neighbourhood
is closed in $X^{ss}$, and is of maximal dimension
equal to the dimension of $G$. 
If $U$ is any $G$-invariant open
subset of the set $X^s$ of stable points of $X$, then $q(U)$ is
an open subset of $X/\!/G$ and the restriction
$q|_U :U \to q(U)$ of $q$ to $U$ is an orbit space
for the action of $G$ on $U$,
so we will write $U/G$ for $q(U)$.
In particular there is an orbit space $X^s/G$ for the action
of $G$ on $X^s$, and $X/\!/G$ is a projective completion of this orbit space.

\begin{equation} \label{gitpicture}
\begin{array}{ccccc}
   X^s & \subseteq & X^{ss} & \subseteq & X \\
       & {\rm open} &  & {\rm open} & \\
  \Big\downarrow & & \Big\downarrow & & \\
  &  &  &  &  \\
  X^s/G & \subseteq &   X/\!/G = X^{ss} / \sim & & \\
       & {\rm open} & & & 
\end{array} \end{equation}

\bigskip

 $X^s, X^{ss}$, and $X/\!/G$ are unaltered if the
line bundle $L$ is replaced by $L^{\otimes k}$ for any $k > 0$
with the induced action of $G$, so it is convenient to allow
 fractional linearisations.

Recall  that  a { categorical quotient} of a variety $X$ under an
action of  $G$ is a $G$-invariant morphism $\phi:X \to Y$ from $X$ to a variety $Y$ 
such that any other $G$-invariant morphism $\tilde{\phi}: X \to
\tilde{Y}$ factors as $\tilde{\phi} = \chi \circ \phi$ for a unique
morphism $\chi:Y \to \tilde{Y}$
\cite[Chapter 2, $\S$4]{New}. An { orbit space} for the action
is a categorical quotient $\phi:X \to Y$ such that each fibre $\phi^{-1}(y)$ is
a single $G$-orbit, and a { geometric
quotient} is an orbit space $\phi:X \to Y$ which is an affine morphism such that

(i) if $U$ is an open affine subset of $Y$ then
$$ \phi^*: {\calo}(U) \to {\calo}(\phi^{-1}(U)) $$
induces an isomorphism of ${\calo}(U)$ onto ${\calo}(\phi^{-1}(U))^G$, and

(ii) if $W_1$ and $W_2$ are disjoint closed $G$-invariant
subvarieties of $X$ then their images $\phi(W_1)$ and $\phi(W_2)$ in
$Y$ are disjoint closed subvarieties of $Y$.

\noindent If $U$ is any
$G$-invariant open subset of the set $X^s= X^s(L)$ of stable points
of $X$, then $q(U)$ is an open subset of $X/\!/G$ and the
restriction $q|_U :U \to q(U)$ of $q$ to $U$ is a
geometric quotient for the action of $G$ on $U$. In particular
$X^s/G = q(X^s)$ is a geometric quotient for the action of $G$ on
$X^s$, while $q: X^{ss} \to X/\!/G$ is a categorical quotient
of $X^{ss}$ under the action of $G$.

The subsets $X^{ss}$ and $X^s$ of $X$ are characterised by the following properties (see
Chapter 2 of \cite{GIT} or \cite{New}). 

\begin{prop} (Hilbert-Mumford criteria)
\label{sss} (i) A point $x \in X$ is semistable (respectively
stable) for the action of $G$ on $X$ if and only if for every
$g\in G$ the point $gx$ is semistable (respectively
stable) for the action of a fixed maximal torus of $G$.

\noindent (ii) A point $x \in X$ with homogeneous coordinates $[x_0:\ldots:x_n]$
in some coordinate system on $\PP^n$
is semistable (respectively stable) for the action of a maximal 
torus of $G$ acting diagonally on $\PP^n$ with
weights $\a_0, \ldots, \a_n$ if and only if the convex hull
$$\conv \{\a_i :x_i \neq 0\}$$
contains $0$ (respectively contains $0$ in its interior).
\end{prop}

In \cite{GIT} the definitions 
of $X^s$ and $X^{ss}$ are extended as follows to allow $L$ to be not ample and
$X$ not projective. However it is not necessarily
the case that $U^{ss} = U \cap X^{ss}$ or that $U^s = U \cap X^s$
when $U$ is a $G$-invariant open subset of $X$.

\begin{defn} \label{defn:s/ssred}
Let $X$ be a { quasi-projective} complex variety with an action of a 
complex reductive
group $G$ and linearisation $L$ on $X$. Then $y \in X$ is {\em
semistable} for this linear action if there exists some $m \geq 0$
and $f \in H^0(X, L^{\otimes m})^G$ not vanishing at $y$ such that
the open subset
$$ X_f := \{ x \in X \ | \ f(x) \neq 0 \}$$
is affine, and $y$ is {\em stable} if also the action of $G$
on $X_f$ is closed with all stabilisers finite.
\end{defn}

When
$X$ is projective and $L$ is ample and $f \in H^0(X, L^{\otimes
m})^G \setminus \{ 0 \}$ for some $m \geq 0$, then $X_f$ is affine
when $f$ is nonconstant, so
this is equivalent to the previous definition.

\begin{rem} \label{reductiveenvquot}
The reason for introducing the requirement that $X_f$ must be affine
in Definition \ref{defn:s/ssred} above is to ensure that $X^{ss}$ has
a categorical quotient $X^{ss} \to \xg$, which restricts to a geometric
quotient $X^s \to X^s/G$ (see \cite{GIT} Theorem 1.10); the quotient $\xg$ 
is quasi-projective, but need not be projective even when $X$ is
projective, if $L$ is not ample. However when $X$ is projective we can define \lq naively
stable' and \lq naively semistable' points by omitting
the condition that $X_f$ should be affine. More precisely,
let $I = \bigcup_{m>0} H^0(X,L^{\otimes m})^G$
and for $f \in I$ let $X_f$ be the $G$-invariant open subset
of $X$ where $f$ does not vanish.
Then a point $x \in X$ is {\em naively semistable} if there exists some 
$f \in I$
which does not vanish at $x$, so that the set of naively semistable points 
is 
$$X^{nss} =  \bigcup_{f \in I} X_f,$$
whereas  $X^{ss} =  \bigcup_{f \in I^{ss}} X_f$ where
$$I^{ss} = \{f
\in I \ | \ X_f
\mbox{ is affine }   \}.$$
The set of {\em naively stable}
points of $X$ is
     $$X^{ns} = \bigcup_{f \in I^{ns}} X_f$$ where
$$I^{ns} = \{f
\in I \ | \  \mbox{ the action of $G$
on $X_f$ is closed with all stabilisers finite} \}$$
while
     $X^{s} = \bigcup_{f \in I^{s}} X_f$ where
$$I^{s} = \{f
\in I^{ss} \ | \  \mbox{ the action of $G$
on $X_f$ is closed with all stabilisers finite} \}.$$
If 
$${\hat{\calo}}_L(X) = \bigoplus_{k \geq 0} H^0(X, L^{\otimes k}) $$
is finitely generated as a complex algebra, then so is its ring of invariants
${\hat{\calo}}_L(X)^G$ because $G$ is reductive, and then 
$$\overline{\xg} = \Proj({\hat{\calo}}_L(X)^G)$$
is a projective completion of $\xg$ and the categorical
quotient $X^{ss} \to \xg$ is the restriction of a
natural $G$-invariant morphism
$$X^{nss} \to \overline{\xg},$$
which by analogy with Definition 3.3 below we will call an
{\em enveloping quotient} of $X^{nss}$.
\end{rem}

Throughout the remainder of $\S$2 we will assume that $X$ is projective and $L$ is ample.

\subsection{Partial desingularisations of quotients}

When $X$ is nonsingular then $X^s/G$ has only orbifold singularities, but if $X^{ss} \neq X^s$
the GIT quotient $\xg$ is likely to have more serious singularities.
However if $X^s \neq \emptyset$ there is a canonical
procedure (see \cite{K2}) for constructing a partial resolution of singularities
$\txg$ of the quotient $\xg$. This involves blowing $X$ up along 
a sequence of nonsingular
$G$-invariant subvarieties, all contained in the complement  of 
the set $X^s$ of stable points of $X$,
to obtain eventually a nonsingular projective variety $\tilde{X}$ with a linear $G$-action,
lifting the action on $X$, for which every semistable point of $\tilde{X}$ is stable.
The blow-down map $\pi:\tilde{X} \to X$ induces a birational morphism $\pi_G:
\txg \to \xg$ which
is an isomorphism over the dense open subset $X^s/G$ of $\xg$,
and if $X$ is nonsingular the quotient $\txg$ has only orbifold singularities.

This construction works as follows \cite{K2}.
Let $V$ be any nonsingular $G$-invariant closed subvariety of $X$ and let
$\pi:\hat{X} \to X$ be the blowup of $X$ along $V$. The linear action of $G$ on the
ample line bundle $L$ over $X$ lifts to a linear action on the line bundle
over $\hat{X}$ which is the pullback of $L^{\otimes k}$ 
tensored with ${\mathcal O}(-E)$,
where $E$ is the exceptional divisor and $k$ is a fixed positive integer.
When $k$ is large the line bundle $\pi^* L^{\otimes k} \otimes {\mathcal O}(-E)$
is ample on $\hat{X}$, and  this linear action satisfies
the following properties:

(i) if $y$ is semistable in $\hat{X}$ then $\pi(y)$ is semistable in $X$;

(ii) if $\pi(y)$ is stable in $X$ then $y$ is stable in $\hat{X}$;

(iii) if $k$ is large enough then the sets $\hat{X}^s$
and $\hat{X}^{ss}$ of stable and semistable points of $\hat{X}$ with respect to this
linearisation are independent of $k$ (cf. \cite{Reichstein}).

\noindent Now $X$ has semistable points which are not stable if and only if there exists
a nontrivial connected reductive subgroup of $G$ which fixes some
semistable point. If so, let $r>0$ be the maximal dimension of the reductive
subgroups of $G$ fixing semistable points of $X$, and let ${\mathcal R}(r)$
be a set of representatives of conjugacy classes in $G$ of all connected 
reductive subgroups $R$ of dimension $r$ such that
$$Z_R^{ss} = \{x\in X^{ss} : \mbox{$R$ fixes $x$} \}$$
is nonempty. Then 
$$\bigcup_{R \in {\mathcal R}(r)} GZ_R^{ss}$$
is a disjoint union of nonsingular closed subvarieties of $X^{ss}$,
and
$$G Z_R^{ss} \cong G \times_{N^R} Z_R^{ss}$$
where $N^R$ is the normaliser of $R$ in $G$.

By Hironaka's theorem we can resolve the singularities of the
closure of $\bigcup_{R \in {\mathcal R}(r)} GZ_R^{ss}$ in $X$ by 
performing a sequence of
blow-ups along nonsingular $G$-invariant closed subvarieties of
$X-X^{ss}$. We then blow up along the proper transform of the
closure of $\bigcup_{R \in {\mathcal R}(r)} GZ_R^{ss}$ to get a nonsingular
projective variety $\hat{X}_1$. The linear action of $G$ on $X$ lifts
to an action on this blow-up $\hat{X}_1$ which can be linearised using
suitable ample line bundles as above,
and it is shown in \cite{K2} that the set $\hat{X}_1^{ss}$ of semistable
points of $\hat{X}_1$ with respect to any of these suitable linearisations
of the lifted action is the complement in the inverse image of $X^{ss}$
of the proper transform of the subset
$$\phi^{-1}\left(\phi\left(\bigcup_{R \in {\mathcal R}(r)} 
GZ_R^{ss}\right)\right)$$
of $X^{ss}$, where $\phi:X^{ss} \to \xg$ is the canonical map. Moreover
no point of $\hat{X}_1^{ss}$ is fixed by a reductive subgroup of $G$ of dimension
at least $r$, and a point in $\hat{X}_1^{ss}$ is fixed by a reductive subgroup
$R$ of $G$ of dimension less than $r$ if and only if it belongs to the
proper transform of the subvariety $Z_R^{ss}$ of $X^{ss}$.

The same procedure can now be applied  to $\hat{X}_1$ to obtain
$\hat{X}_2$ such that no reductive subgroup of $G$ of dimension at least $r-1$
fixes a point of $\hat{X}_2^{ss}$. After repeating enough times we obtain 
$\tilde{X}$ satisfying $\tilde{X}^{ss}=\tilde{X}^s$.
 If we are only interested in $\tilde{X}^{ss}$ and the partial
resolution $\txg$ of $\xg$, rather than in $\tilde{X}$ itself, then there
is no need in this procedure to resolve the singularities of the closure
of $\bigcup_{R \in {\mathcal R}(r)} GZ_R^{ss}$ in $X$. Instead we can simply blow
$X^{ss}$ up along $\bigcup_{R \in {\mathcal R}(r)} GZ_R^{ss}$ (or equivalently
along each $GZ_R^{ss}$ in turn) and let $\hat{X}_1^{ss}$ be the set of semistable
points in the result, and then repeat the process.

Thus the geometric quotient $X^s/G$ has {\em two} natural compactifications
$\xg$ and $\widetilde{\xg} = \txg$, which fit into a diagram 

$$ \begin{array}{cccccccccc}
X^{{s}}/G & \hookrightarrow & \widetilde{\xg} = \tilde{X}^{ss}/G &
\leftarrow & \tilde{X}^{ss} = \tilde{X}^s & \hookrightarrow & \tilde{X}\\
|| & &  \downarrow &  & \downarrow & &
\downarrow \\ X^s/G & \hookrightarrow & \xg 
& \leftarrow & X^{ss} & \hookrightarrow & X.
\end{array}
$$

\subsection{Variation of GIT}

The GIT quotient $X/\!/G$ depends not just on the action of
$G$ on $X$ but also on the choice of linearisation $\call$ of the
action; that is, the choice of the line bundle $L$ and
the lift of the action to $L$. This should of course
be reflected in the notation; to avoid ambiguity we will
sometimes add appropriate decorations, as in
$X/\!/_{\call}G$ and also $X^{s,\call}$ and
$X^{ss,\call}$. There is thus a natural question:
how does the GIT quotient $X/\!/_{\call} G$ vary as the linearisation $\call$
varies? This has been studied by Brion and Procesi \cite{BP}
and Goresky and MacPherson \cite{GM}
(in the abelian case) and by Thaddeus \cite{Thad}, 
Dolgachev and Hu \cite{DolgHu} and Ressayre \cite{Ress}
for general reductive groups $G$.

A very simple case is when the line bundle $L$ is fixed, but
the lift $\tau:G \times L \to L$ of the action of $G$
on $X$ to an action on $L$ varies. The only possible such
variation is to replace $\tau$ with
$$(g,\ell) \mapsto \chi(g) \tau(g,\ell)$$
where $\chi:G \to \CC^*$ is a character of $G$. The 
Hilbert-Mumford criteria (Proposition \ref{sss} above)
can be used to see how $X^s$ and $X^{ss}$ are affected
by such a variation.

\begin{ex}
Consider the linear action of $\CC^*$ on $X=\PP^{r+s}$, with respect to
the hyperplane line bundle $L$ on $X$, where the linearisation $\call_+$
is given by the representation
$$t \mapsto \diag(t^3,\ldots,t^3,t,t^{-1},\ldots,t^{-1})$$
of $\CC^*$ in $GL(r+s+1;\CC)$ in which $t^3$ occurs with
multiplicity $r$ and $t^{-1}$ occurs $s$ times. 
For this linearisation we have
$$X^{ss,\call_+} = X^{s,\call_+} = \{ [x_0:\ldots:x_{r+s}] \in \PP^{r+s}: \mbox{$x_0,\ldots,x_r$
 are not all 0 $$
 $$ \hfill and $x_{r+1},\ldots,x_{r+s}$ are not all 0} \}$$
and $X/\!/_{\call_+} \CC^*$ is isomorphic to the product of a weighted
projective space $\PP(3,\ldots,3,1)$ of dimension $r$ and
the projective space $\PP^{s-1}$.
The same action of $\CC^*$ on $X$
has other linearisations with respect to the hyperplane line bundle;
let $\call_0$ and $\call_-$ denote the linearisations given by 
multiplying the representation above by the characters
$\chi(t) = t^{-1}$ and  $\chi(t) = t^{-2}$.
Then $\call_-$ is given by the representation
$$t \mapsto \diag(t,\ldots,t,t^{-1},t^{-3},\ldots,t^{-3})$$
of $\CC^*$ in $GL(r+s+1;\CC)$ and 
for this linearisation
$$X^{ss,\call_-} = X^{s,\call_-} = \{ [x_0:\ldots:x_{r+s}] \in \PP^{r+s}: \mbox{$x_0,\ldots,x_{r-1}$
 are not all 0 $$
 $$ \hfill and $x_{r},\ldots,x_{r+s}$ are not all 0} \}$$
while $X/\!/_{\call_-} \CC^*$ is isomorphic to the product of 
the projective space $\PP^{r-1}$ and a weighted
projective space $\PP(1,3,\ldots,3)$ of dimension $s$.
Finally for the linearisation $\call_0$ given by the representation
$$t \mapsto \diag(t^2,\ldots,t^2,1,t^{-2},\ldots,t^{-2})$$
of $\CC^*$ in $GL(r+s+1;\CC)$ we have
$$X^{ss,\call_0} \neq X^{s,\call_0}$$
(semistability does not imply stability) and the quotient
$X/\!/_{\call_0} \CC^*$ has more serious singularities than
the orbifolds $X/\!/_{\call_\pm} \CC^*$. It can be identified with the result
of collapsing $[0:\ldots :0:1] \times \PP^{s-1}$ to
a point in $\PP(3,\ldots,3,1) \times \PP^{s-1}$, and also
with the result of collapsing $\PP^{r-1} \times [1:0:\ldots :0]$ to
a point in $\PP^{r-1} \times \PP(1,3,\ldots,3)$.
\label{exab}
\end{ex}

\begin{rem} \label{flips}
The general case when the line bundle $L$ is allowed to vary, as
well as the lift of the $G$-action from $X$ to the line bundle,
can be reduced to the case above by a trick due to Thaddeus \cite{Thad}.
Suppose that a given action of $G$ on $X$ lifts to ample line
bundles $L_0, \dots, L_m$ giving linearisations $\call_0, \ldots, \call_m$ over $X$, and consider
the projective variety
$$Y = \PP(L_0 \oplus \cdots \oplus L_m).$$
Then the induced action of $G$ on $Y$ has a natural linearisation,
and the complex torus $T=\CC^{m+1}$ also acts on $Y$, commuting
with the action of $G$, with a natural linearisation which
can be modified using any character $\chi$ of $T$. Taking
$\chi$ to be the $j$th projection $\chi_j:T \to \CC^*$
and using the fact that the GIT quotient operations with
respect to $G$ and $T$ commute, we find that
$$X/\!/_{\call_j} G \cong (Y/\!/G)/\!/_{\chi_j} T$$
and hence the general question of variation of GIT
quotients with linearisations reduces to the special case
when the variation is by multiplication by a character
of the group. The conclusion \cite{DolgHu,Ress,Thad} is that, roughly speaking,
the space of all possible ample fractional linearisations 
of a given $G$-action on a projective variety $X$ is divided into finitely many polyhedral
chambers within which the GIT quotient is constant. Moreover,
when a wall between two chambers is crossed, the quotient undergoes
a transformation which typically
can be thought of as a blow-up followed by
a blow-down. 
If $\call_+$ and $\call_-$ represent fractional linearisations
in the interiors of two adjoining chambers, and $\call_0$ represents
a fractional linearisation in the interior of the wall between them, then
we have inclusions
$$X^{s,\call_0} \subseteq X^{s,\call_+} \cap X^{s,\call_-}
\   \mbox{ and } \  X^{ss,\call_+} \cup X^{ss,\call_-} \subseteq
X^{ss,\call_0}$$
inducing morphisms
$$X/\!/_{\call_\pm} G \to X/\!/_{\call_0} G$$
which are isomorphisms over $X^{s,\call_0}/G$.
In addition, under mild conditions there are sheaves of ideals on the two
quotients $X/\!/_{\call_\pm} G$ whose blow-ups are both isomorphic to a component
of the fibred product of the two quotients over the
quotient $X/\!/_{\call_0} G$ on the wall. \end{rem}

\begin{rem}
If $\phi:X \to Y$ is a categorical quotient for a $G$-action on a 
variety $X$ then its restriction to a $G$-invariant open subset of $X$ is
not necessarily a categorical quotient for the action of $G$ on $U$.
In the situation above, 
as in Example \ref{exab}, we have $X^{ss,\call_+} \subseteq
X^{ss,\call_0}$ but the restriction of the categorical
quotient $X^{ss,\call_0} \to X/\!/_{\call_0} G$
to $X^{ss,\call_+}$ is not a categorical quotient for 
the $G$ action on $X^{ss,\call_+}$.
\end{rem}

\section{Quotients by non-reductive actions}

Translation actions appear all over geometry, so it is not surprising that
there are many cases
of moduli problems which involve non-reductive group actions,
where Mumford's GIT does not apply. One example
is that of  hypersurfaces in a toric variety $Y$. The case we shall
consider in detail in this paper is when $Y$ is
 the weighted projective plane $\PP(1,1,2)$ (cf. \cite{F3}), with
homogeneous coordinates $x,y,z$ (that is, 
$Y$ is the quotient of $\CC^3 \setminus \{ 0 \}$ by the
action of $\CC^*$ with weights 1,1 and 2, and $x,y$ and $z$ are coordinates
on $\CC^3 \setminus \{ 0 \}$).
Let $H$ be  the
automorphism group of $Y=\PP(1,1,2)$, which is the quotient
by $\CC^*$ of a semidirect product of
the unipotent group $U=(\CC^+)^3$ acting on $Y$ via
$$[x:y:z] \mapsto [x:y:z + \lambda x^2 + \mu xy + \nu y^2] \quad \mbox{
for } (\lambda,\mu,\nu) \in (\CC^+)^3$$ and the reductive group
$GL(2;\CC) \times GL(1; \mathbb{C})$ acting on the $(x,y)$
coordinates and the $z$ coordinate.  
$H$ acts linearly on the projective space $X_d$ of weighted degree $d$
polynomials in $x,y,z$.

\begin{ex} When $d=4$, a basis for the weighted degree $d$
polynomials is
$$\{ x^4, x^3y, x^2y^2,xy^3,y^4,x^2z,xyz,y^2z,z^2 \},$$
and with respect to this basis, the $U$-action is given by
$$\left( \begin{array}{ccccccccc}

1 & 0 & 0 & 0 & 0 & \lambda & 0 & 0 & \lambda^2 \\

0 & 1 & 0 & 0 & 0 & \mu & \lambda & 0 & 2\lambda\mu \\

0 & 0 & 1 & 0 & 0 & \nu & \mu & \lambda & 2\lambda\nu + \mu^2 \\

0 & 0 & 0 & 1 & 0 & 0 & \nu & \mu & 2\mu\nu \\

0 & 0 & 0 & 0 & 1 & 0 & 0 & \nu & \nu^2 \\

0 & 0 & 0 & 0 & 0 & 1 & 0 & 0 & 2\lambda \\

0 & 0 & 0 & 0 & 0 & 0 & 1 & 0 & 2\mu \\

0 & 0 & 0 & 0 & 0 & 0 & 0 & 1 & 2\nu \\

0 & 0 & 0 & 0 & 0 & 0 & 0 & 0 & 1

\end{array}  \right).$$
\end{ex}

The tautological family $\calh^{(d)}$ parametrised 
by $X_d$ of hypersurfaces in $Y$ has
 the
following two properties:

(i) the hypersurfaces ${\cal H}^{(d)}_s$ and ${\cal H}^{(d)}_t$ parametrised
by weighted degree $d$ polynomials $s$ and $t$  are isomorphic 
as hypersurfaces in $Y$ if and
only if $s$ and $t$ lie in the same orbit of the natural action
of $H \cong U \rtimes GL(2;\CC)$ on $X_d$, and

(ii) (local universal property) any family of hypersurfaces in
$Y$ is locally equivalent to the pullback of ${\cal H}^{(d)}$ along
a morphism to $X_d$.

\noindent This means that the construction of a (coarse) moduli space 
of weighted degree $d$ hypersurfaces in $Y$ is equivalent
to constructing an orbit space  for the action of $H$ on
$X_d$ (\cite{New} Proposition 2.13). 

Now let $H$ be any affine algebraic 
group, with unipotent radical $U$, acting linearly on a complex projective variety $X$ with respect to an ample line bundle $L$. Of course the most immediate
difficulty when trying to generalise Mumford's GIT to a non-reductive situation
is that the ring of invariants 
$${\hat{\calo}}_L(X)^H = \bigoplus_{k \geq 0} H^0(X, L^{\otimes k})^H$$
is not necessarily finitely generated as a graded complex algebra,
so that $\Proj({\hat{\calo}}_L(X)^H)$ is not well-defined as a projective variety.
However $\Proj({\hat{\calo}}_L(X)^H)$ does make sense as a scheme, and the
inclusion of ${\hat{\calo}}_L(X)^H$ in ${\hat{\calo}}_L(X)$ gives us a rational map of schemes 
$q$ from $ X$ to $ \Proj({\hat{\calo}}_L(X)^H)$, whose
image is a constructible subset of $\Proj({\hat{\calo}}_L(X)^H)$ (that is, a finite union of
locally closed subschemes). 

The action on $X$ of the unipotent radical $U$ of $H$
is studied in \cite{DK} (building on earlier work such as \cite{F2,F1,GP1,GP2,W}), where the following definitions
are made and results proved.

\begin{defn}
Let $I = \bigcup_{m>0} H^0(X,L^{\otimes m})^U$
and for $f \in I$ let $X_f$ be the $U$-invariant affine open subset
of $X$ where $f$ does not vanish, with ${\calo}(X_f)$ its coordinate ring.
A point $x \in X$ is {\em naively semistable} if the rational map $q$ from
$X$ to $ \Proj({\hat{\calo}}_L(X)^U)$ is
well-defined at $x$; that is, if there exists some $f \in I$
which does not vanish at $x$.  The set of naively semistable points 
is $X^{nss}= \bigcup_{f \in I} X_f$.
The {\em finitely generated semistable set} of $X$ is  $$X^{ss,
fg} =  \bigcup_{f \in I^{fg}} X_f$$ where
$$I^{fg} = \{f
\in I \ | \ {\calo}(X_f)^U
\mbox{ is finitely generated }   \}.$$
The set of {\em naively stable}
points of $X$ is
     $$X^{ns} = \bigcup_{f \in I^{ns}} X_f$$ where
$$I^{ns} = \{f
\in I \ | \ {\calo}(X_f)^U
\mbox{ is finitely generated, and }$$
$$  q: X_f \longrightarrow
\Spec({\calo}(X_f)^U) \mbox{ is a geometric quotient} \}.$$
The set of {\em locally trivial stable} points is $$ X^{lts} =
\bigcup_{f \in I^{lts} } X_f$$ where
\begin{eqnarray*} I^{lts}\ \  =\ \  \{f
\in I  \ | \  {\calo}(X_f)^U
\mbox{
is finitely generated, and \ \ \ \ \ \ } \\
    q: X_f \longrightarrow \Spec({\calo}(X_f)^U) \mbox{ is a locally trivial
geometric quotient} \}. \end{eqnarray*}

\end{defn}

\begin{defn}\label{defn:envelopquot}
Let $q: X^{ss, fg} \rightarrow \Proj({\hat{\calo}}_L(X)^U)$ be the natural
morphism of schemes.  The {\em enveloped quotient} of
$X^{ss,fg}$ is $q: X^{ss, fg} \rightarrow q(X^{ss,fg})$, where
$q(X^{ss,fg})$ is a dense constructible subset of the {\em
enveloping quotient}
$$X /\!/ U = \bigcup_{f \in I^{ss,fg}}
\Spec({\calo}(X_f)^U)$$ of $X^{ss, fg}$. Note that $q(X^{ss,fg})$ is 
not necessarily a subvariety of $X /\!/ U $, as is demonstrated by
the example studied in \cite{DK} $\S$6 of $U = \CC^+$ acting on $X=\PP^n$
via the $n$th symmetric product of its standard representation on $\CC^2$
when $n$ is even.
\end{defn}

\begin{prop} {\rm (\cite{DK} 4.2.9 and 4.2.10).} \label{patching2} The enveloping quotient
$X/\!/U$ is a quasi-projective variety with an ample line bundle
$L_H \to X/\!/U$ which pulls back to a positive tensor power of $L$
under the natural map $q:X^{ss,fg} \to X/\!/U$.
If ${\hat{\calo}}_L(X)^U$ is finitely generated then $X/\!/U$ is the projective
variety $\Proj({\hat{\calo}}_L(X)^U)$.
\end{prop}

Now suppose that $G$ is a complex reductive group with  $U$
as a closed subgroup.   Let $G \times_U X$ denote the quotient of $G \times X$
by the free action of $U$ defined by $u(g,x)=(g u^{-1}, ux)$ for $u \in U$,
which is a quasi-projective variety by \cite{PopVin} Theorem 4.19. There
is an induced $G$-action on $G \times_U X$ given by left
multiplication of $G$ on itself.
If the action of $U$ on $X$ extends to an action of $G$ there is an isomorphism of
$G$-varieties 
$$G \times_U X \cong (G/U) \times X$$ given by
\begin{equation} \label{9Febiso} [g,x] \mapsto (gU, gx). \end{equation}

When $U$ acts linearly on $X$ with respect to a very ample line bundle $L$ inducing an
embedding of $X$ in $\PP^n$, and $G$ is a subgroup of $SL(n+1; \CC)$,
then we get  a very ample $G$-linearisation 
(which by abuse of notation we will also denote by $L$) on $G \times_U X$
as follows:
$$ G \times_U X \hookrightarrow G \times_U \mathbb{P}^n  \cong
(G/U) \times \mathbb{P}^n,
$$by  taking the trivial bundle on the quasi-affine variety
$G/U$. If we choose a $G$-equivariant embedding of $G/U$ in an affine
space $\mathbb{A}^m$ with a linear $G$-action we get a $G$-equivariant
embedding of $G \times_U X$ in
$$\mathbb{A}^m \times \mathbb{P}^n \subset \mathbb{P}^m \times \mathbb{P}^n
\subset \mathbb{P}^{nm+m+n}$$
and the $G$-invariants on $ G \times_U X$ are given by
\begin{equation} \label{name}  \bigoplus_{m \geq 0} H^0(G \times_U X, L^{\otimes m})^G \cong
\bigoplus_{m \geq 0} H^0(X, L^{\otimes m})^U = {\hat{\calo}}_L(X)^U.\end{equation}

\begin{defn}
The sets of {\em Mumford stable points} and
{\em Mumford semistable points} in $X$ are $X^{ms} = i^{-1}((G \times_U
X)^s)$ and $X^{mss} = i^{-1}((G \times_U
X)^{ss})$
where $i: X \rightarrow G \times_U X$ is the inclusion given by
$x \mapsto [e,x]$ for $e$ the identity element of $G$.
 Here $(G \times_U X)^s$ and $(G \times_U X)^{ss}$ are
defined as in Definition \ref{defn:s/ssred} for the induced linear
action of $G$ on the quasi-projective variety $G \times_U X$.
\end{defn}

In fact it follows from Theorem \ref{thm:main} below that
$X^{ms}$ and $X^{mss}$ are equal and are independent of the choice of $G$.

\begin{defn} \label{defn:separ}
 A {\em finite separating
set of invariants} for the linear action
of $U$ on $X$ is a collection of invariant sections $\{f_1,
\ldots, f_n \}$ of positive tensor powers of $L$ such that, if $x,y$
are any two points of $X$ then $ f(x) = f(y)$ for all
invariant sections $f$
of $L^{\otimes k}$ and all $k 
> 0$ if and only if
$$ f_i(x) =
f_i(y) \hspace{.5in} \forall i = 1, \ldots, n.
$$
If $G$ is any reductive group containing
$U$, a finite separating set $S$ 
of invariant sections of positive tensor powers of $L$ is a {\em
finite fully separating set of invariants} for the linear $U$-action
on $X$ if

(i) for every $x \in X^{ms}$ there exists $f \in S$ with associated
$G$-invariant $F$ over $G \times_U X$ (under the isomorphism
(\ref{name})) such that $x \in (G \times_U
X)_{F}$ and $(G \times_U X)_{F}$ is affine; and

(ii) for every $x \in X^{ss,fg}$ there exists $f \in S$ such that $x
\in X_f$ and $S$ is a generating set for ${\calo}(X_f)^U$.
\end{defn}

This
definition is in fact independent of the choice of $G$ (see \cite{DK} Remark 5.2.3).

\begin{defn}\label{defn:envelope}
Let $X$ be a quasi-projective variety with a linear $U$-action with
respect to an ample line bundle $L$ on $X$, and let $G$ be a complex reductive group
containing $U$ as a closed subgroup. A
$G$-equivariant projective completion $\overline{G \times_U X}$ of
$G \times_U X$, together with a $G$-linearisation with respect to
a line bundle $L$ which
restricts to the given $U$-linearisation on $X$,
is a
{\em reductive envelope} of the linear $U$-action on $X$ if every
$U$-invariant $f$ in some finite fully separating set of invariants
$S$ for the $U$-action on $X$  extends to a $G$-invariant section of a tensor
power of $L$ over $\overline{G \times_U X}$.

If moreover there exists such an $S$ for which every
$f\in S$ extends to a $G$-invariant section $F$ over $\overline{G
\times_U X}$ such that $(\overline{G \times_U X})_F$ is affine, then
 $(\overline{G \times_U X}, L')$ is a {\em fine reductive
envelope}, and if $L$ is ample (in which case
$(\overline{G \times_U X})_F$ is always affine) it is an {\em ample
reductive envelope}.

If every $f \in S$
extends to a $G$-invariant $F$ over $\overline{G \times_U X}$ which
vanishes on each 
 codimension 1 component of the boundary of $G
\times_U X$ in $\overline{G \times_U X}$, then 
a reductive envelope for the linear
$U$-action on $X$ is called a {\em strong} reductive envelope.
\end{defn}

It will be useful to add an extra definition which does not
appear in \cite{DK}.

\begin{defn} \label{defn:stablyfine}
In the notation of Definitions \ref{defn:separ} and
\ref{defn:envelope} above, a reductive envelope is called {\em stably fine} if for
every $x \in X^{ms}$ there exists a $U$-invariant $f$ which extends
to a $G$-invariant section $F$ over $\overline{G \times_U X}$ such that 
$x \in (G \times_U X)_F$ and both $(G \times_U X)_F$ and $(\overline{G \times_U X})_F$
are affine.
\end{defn}

\begin{defn}\label{defn:s/ssbar} Let $X$ be a
projective variety with a linear $U$-action and a reductive envelope
$\overline{G \times_U X}$.  The set of {\em completely stable points} of $X$ with
respect to the reductive envelope is $$X^{\overline{s}} = (j
\circ i)^{-1}(\overline{G \times_U X}^s)$$ and the set of {\em
completely semistable points} is $$X^{\overline{ss}} = (j \circ
i)^{-1}(\overline{G \times_U X}^{ss}),$$
where
 $i: X \hookrightarrow G
\times_U X$ and $j: G \times_U X \hookrightarrow \overline{G
\times_U X}$ are the inclusions, and $\overline{G \times_U X}^{s}$
and $\overline{G \times_U X}^{ss}$ are the stable and semistable sets for
the linear $G$-action on $\overline{G \times_U X}$. Following Remark
\ref{reductiveenvquot} we also define
$$X^{\overline{nss}} = (j \circ
i)^{-1}(\overline{G \times_U X}^{nss});$$
then $X^{\overline{nss}} = X^{\overline{ss}}$ when the reductive envelope
is ample, but not in general otherwise.
\end{defn}

\begin{theorem}\label{thm:main} {\rm (\cite{DK} 5.3.1 and 5.3.5).}
Let $X$ be a normal projective variety with a linear $U$-action, for $U$ a
connected unipotent group, and let $(\overline{G \times_U X},L)$ be
any fine reductive envelope. Then
$$
X^{\overline{s}}  \subseteq  X^{lts} = X^{ms} = X^{mss}
\subseteq  X^{ns}   \subseteq  X^{ss, fg}  \subseteq X^{\overline{ss}} = X^{nss}.
$$
The stable sets $X^{\overline{s}}$,
$X^{lts} = X^{ms} = X^{mss}$ and $X^{ns}$ admit quasi-projective
geometric quotients, given by restrictions of the quotient map $q =
\pi \circ j \circ i$
where
$$\pi: (\overline{G \times_U X})^{ss} \to
\overline{G \times_U X}/\!/G$$
is the classical GIT quotient map for
the reductive envelope and $i,j$ are as in Definition \ref{defn:s/ssbar}.
 The quotient map $q $
restricted to the open subvariety $X^{ss,fg}$ is an enveloped
quotient with $q: X^{ss,fg} \rightarrow X /\!/ U$ an enveloping
quotient. Moreover
  $X /\!/ U$ is an open
subvariety of $\overline{G \times_U X}/\!/G$ and there is an ample
line bundle $L_U$ on $X/\!/U$ which pulls back to a tensor power
$L^{\otimes k}$ of the line bundle $L$ for some $k>0$
and extends to an ample line bundle on 
$\overline{G \times_U X}/\!/G$.

If furthermore $\overline{G \times_U X}$ is normal and  provides a fine strong reductive envelope for the linear
$U$-action on $X$, then $X^{\overline{s}} = X^{lts}$ and $X^{ss,fg} =
X^{nss}$.
\end{theorem}

\begin{defn}\label{defn:main} {\rm (\cite{DK} 5.3.7).}
Let $X$ be a projective variety
equipped with a linear $U$-action. A point $x \in X$ is called
{\em stable} for the linear $U$-action if $x \in X^{lts}$ and {\em semistable} if $x
\in X^{ss, fg}$, so from now on we will write $X^s$ (or $X^{s,U}$) for $X^{lts}$
and $X^{ss}$ (or $X^{ss,U}$) for $X^{ss,fg}$.
\end{defn}

Thus in the situation of Theorem \ref{thm:main} we have a
diagram of quasi-projective varieties
$$ \begin{array}{ccccccccccccc}
X^{\overline{s}} & \subseteq & X^{s} &
\subseteq & X^{ns} & \subseteq & X^{ss} & \subseteq & X^{\overline{ss}} = X^{nss}\\
\downarrow &  & \downarrow &  & \downarrow & &
\downarrow & & \downarrow \\ X^{\overline{s}}/U & \subseteq &
X^{s}/U & \subseteq & X^{ns}/U & \subseteq & X/\!/U & \subseteq &
\overline{G \times_U X}/\!/G
\end{array}
$$
where all the inclusions are open and all the
vertical morphisms are restrictions of
$\pi:(\overline{G \times_U X})^{ss} \to \overline{G \times_U X}/\!/G$,
and each except the last is a restriction of
the map of schemes $q:X^{nss} \to \Proj({\hat{\calo}}_L(X)^U)$) associated to
the inclusion ${\hat{\calo}}_L(X)^U \subseteq {\hat{\calo}}_L(X)$.
In particular we have

\begin{equation} \label{xudiagram} \begin{array}{ccccccccccccc}
 X^{s} &
\subseteq &  X^{ss} & \subseteq & X \\
\downarrow &  & \downarrow & & \\ 
X^{s}/U & \subseteq &  X/\!/U & & 
\end{array}
\end{equation}
which looks very similar to the situation for reductive actions
(see diagram (\ref{gitpicture}) above), with
the major differences that

(i) $\xu$ is not always projective,

\noindent and (even if the ring of invariants ${\hat{\calo}}_L(X)^U$ is finitely generated
and $\xu = \Proj ({\hat{\calo}}_L(X)^U)$ is projective)

(ii) the morphism $X^{ss} \to \xu$ is not in general surjective.

\begin{rem} \label{rem:stablyfine}
The proofs of \cite{DK} Theorem 5.3.1 and Theorem 5.3.5 show that if $X$ is
a normal projective variety with a linear $U$-action and
$(\overline{G \times_U X},L)$ is any stably fine
reductive envelope in the sense of Definition \ref{defn:stablyfine},
then 
$$X^{\bar{s}} \subseteq X^s \subseteq X^{ns} \subseteq X^{ss}$$
and if furthermore $\overline{G \times_U X}$ is normal and
provides a reductive envelope which is both strong and
stably fine, then $X^{\bar{s}} = X^s$. Indeed, this is still true even if
$(\overline{G \times_U X},L)$ is not a reductive envelope at all,
provided that it satisfies all the conditions except for omitting (ii) 
in Definition \ref{defn:separ}.
\end{rem}

There always exists an ample, and hence fine, but not necessarily strong, reductive envelope 
for any linear $U$-action on a projective variety $X$, at least if we replace the line bundle $L$ with a suitable
positive tensor power of itself, by \cite{DK} Proposition
5.2.8. By Theorem \ref{thm:main} above a choice of
fine reductive envelope $\overline{G \times_U X}$ provides a
projective completion
$$\overline{X/\!/U} = \overline{G \times_U X}/\!/G$$
of the enveloping quotient $X/\!/U$. This projective
completion in general depends on the
choice of reductive envelope, but when ${\hat{\calo}}_L(X)^U$ is
finitely generated then $X/\!/U = \Proj({\hat{\calo}}_L(X)^U)$ is itself
projective, which implies that
${X/\!/U} = \overline{G \times_U X}/\!/G$
for any fine reductive envelope $\overline{G \times_U X}$.

The proof of \cite{DK} Theorem 5.3.1 also gives us the following result for
any reductive envelope, not necessarily fine or strong.

\begin{prop}\label{prop:main} 
Let $X$ be a normal projective variety with a linear $U$-action, for $U$ a
connected unipotent group, and let $(\overline{G \times_U X},L)$ be
any reductive envelope. Then
$$
X^{\overline{s}} \subseteq X^s \subseteq X^{ss}  \subseteq X^{\overline{nss}}, 
$$
and if the graded algebra $\bigoplus_{k \geq 0} H^0(\overline{G \times_U X},L^{\otimes k})$ is finitely
generated then the projective completion
$$\overline{\overline{G \times_U X}/\!/G} = \Proj( \bigoplus_{k \geq 0} H^0(\overline{G \times_U X},L^{\otimes k}) )$$ 
of  $\overline{G \times_U X}/\!/G$   (cf. Remark \ref{reductiveenvquot})
is a projective completion of $\xu$ with a commutative
diagram 
\begin{equation} \label{xudiagram2} \begin{array}{ccccccccccccc}
X^{\overline{s}} & \subseteq & X^s & \subseteq & X^{ss} &
\subseteq &  X^{\overline{nss}}  \\
\downarrow &  & \downarrow  & & \downarrow & & \downarrow \\ X^{\overline{s}}/U & \subseteq & X^s/U & \subseteq &
 X/\!/U & \subseteq  & \overline{X/\!/U} = \overline{\overline{G \times_U X}/\!/G}.
\end{array}
\end{equation}
\end{prop}

\section{Choosing reductive envelopes}

Let $H$ be a connected affine algebraic group over $\CC$. Then $H$ has a unipotent radical $U$, which
is a normal subgroup of $H$ with reductive quotient group $R=H/U$.
We can hope to quotient first by the action of $U$, and then
by the induced action of the reductive group $H/U$, provided that
the unipotent quotient is sufficiently canonical to inherit an
induced linear action of the reductive group $R$. Moreover $U$
has canonical series of normal subgroups $\{1\}=U_0 \leq U_1 \leq \cdots
\leq U_s=U$ such that each
successive subquotient is isomorphic to $(\CC^+)^r$ for some $r$
(for example the descending central series of $U$),
so we can hope to quotient successively by unipotent
groups of the form $(\CC^+)^r$, and then finally by the reductive
group $R$. Therefore we will concentrate on the case
when $U \cong (\CC^+)^r$ for some $r$; of course this is
the situation in our example concerning hypersurfaces in the
weighted projective plane $\PP(1,1,2)$, when $H$ is the automorphism
group of $\PP(1,1,2)$ and $U$ is its unipotent radical. 

More generally, let us assume first that $U$ is a
unipotent group with a one-parameter group of automorphisms
$\lambda:\CC^* \to \mbox{Aut}(U)$ such that the weights of
the induced $\CC^*$ action on the Lie algebra $\lieu$ of
$U$ are all nonzero.  
When $U = \cplusr$ we can
take $\lambda$ to be the inclusion of the central
$\CC^*$ in $\mbox{Aut}(U) \cong \glr$. Then we can form the semidirect
product
$$\hat{U} = \CC^* \ltimes U$$
given by $\CC^* \times U$ with group multiplication
$$(z_1,u_1).(z_2,u_2) = (z_1 z_2, (\lambda(z_2^{-1})(u_1))u_2).$$ 
$U$ meets the centre
of $\hat{U}$ trivially, so we have an inclusion
$$ U \hookrightarrow \hat{U} \to \mbox{Aut}(\hat{U})
\to GL(\mbox{Lie}\hat{U}) = GL(\CC \oplus \lieu)$$
where $\hat{U}$ maps to its group of inner automorphisms.
Thus $U$ is isomorphic to a closed subgroup of the reductive
group $G=SL(\CC \oplus \lieu)$.

In particular when $U = \cplusr$ we have $U \leq G = \slrplus$, and then
$$G/U \cong \{ \alpha \in (\CC^r)^* \otimes \CC^{r+1} |  \ \alpha:\CC^r \to \CC^{r+1} 
\mbox{ is injective }\}$$
with the natural $G$-action $g\alpha = g \circ \alpha $. Since the injective
linear maps from $\CC^r$ to $\CC^{r+1}$ form an open subset in the affine
space $(\CC^r)^* \otimes \CC^{r+1}$ whose complement has codimension two,
it follows that $U = \cplusr$ is a Grosshans subgroup of $G = \slrplus$
and 
$$\calo(G)^U \cong \calo(G/U) \cong \calo((\CC^r)^* \otimes \CC^{r+1})$$
is finitely generated \cite{Grosshans}. 

\subsection{Actions of $\cplusr$ which extend to $\slrplus$}
Let $X$ be a normal projective variety with a linear action
of $U = \cplusr$ with respect to an ample line bundle $L$.
Suppose first that the linear action of $U = \cplusr$
on $X$ extends to a linear action of $G = \slrplus$, giving us an identification of
$G$-spaces
$$G \times_U X \cong (G/U) \times X$$ as at
(\ref{9Febiso}) via $ [g,x] \mapsto (gU, gx)$. Then 
(as in the Borel transfer theorem \cite[Lemma 4.1]{Dolg})
$$\hat{\calo}_L(X)^U \cong \hat{\calo}_L(G \times_U X)^G \cong [\calo(G/U) \otimes \hat{\calo}_L(X)]^G$$
is finitely generated \cite{Grosshans2} and we have a
reductive envelope
$$\overline{G \times_U X} = \prrplus
\times X,$$
where $\prrplus = {\PP(\CC \oplus ((\CC^r)^* \otimes \CC^{r+1}))}$,
with 
$$\overline{G \times_U X}/\!/G \cong X/\!/U = \Proj(\hat{\calo}_L(X)^U).$$
More precisely, if we choose  for our linearisation on $\overline{G \times_U X}$ the line bundle
$$ L^{(N)} = \calo_\prrplus (N) \otimes L$$
with $N>0$ sufficiently large,
then by \cite{DK} Lemma 5.3.14 we obtain a reductive envelope which is strong as well
as ample, and so by Theorem \ref{thm:main} we have
\begin{equation}
X^{\bar{s}} = X^s \mbox{ and } X^{\bar{ss}} = X^{ss}. \label{sbarssbar} \end{equation}

\begin{rem}
Even if $X$ is nonsingular, this quotient 
$$\xu = (\prrplus \times X) /\!/G$$
may have serious singularities if there are semistable points which
are not stable. However provided that $X^s \neq \emptyset$ we can construct a partial
desingularisation
$${\widetilde{\xu}}^{(G)} = (\widetilde{\prrplus \times X}) /\!/G$$
as in $\S$2.2 by blowing 
$\prrplus \times X $  up successively along $G$-invariant closed subvarieties,
all disjoint from $(\prrplus \times X)^s$ and hence from $X^{\overline{s}} = X^s$, to get a
linear $G$-action on the resulting blow-up
$\widetilde{\prrplus \times X}$ 
for which all semistable points are stable.
This construction is determined by the linear $G$-action, and if $X$ is
nonsingular the resulting quotient 
is an orbifold. 
Since $X^{\bar{s}} = X^s$ and the morphism
$${\widetilde{\xu}}^{(G)} = (\widetilde{\prrplus \times X}) /\!/G \to 
(\prrplus \times X)/\!/G = \xu$$
is an isomorphism over $(\prrplus \times X)^s$, it follows
that ${\widetilde{\xu}}^{(G)} \to \xu$ is an isomorphism over $X^s/U$, and
hence we have two compactifications of the geometric quotient $X^s/U$:
$$ \begin{array}{ccccccccc}
X^{{s}}/U & \subseteq & {\widetilde{\xu}}^{(G)}\\
|| &  & \downarrow \\ X^{{s}}/U & \subseteq &
\xu
\end{array}
$$
where ${\widetilde{\xu}}^{(G)}$ is an orbifold.
\end{rem}

$U = \cplusr$ is the unipotent radical of a parabolic
subgroup 
\begin{equation} \label{parabolic} P = U   \rtimes  \glr \end{equation}
 in $\slrplus$ with Levi subgroup
$\glr$ embedded in $\slrplus$ as
$$ g \mapsto \left( \begin{array}{ccc}
g & 0 \\ 
0 & \det g^{-1}
\end{array}\right) .$$
We have
$$ G \times_U X \cong G \times_P (P \times_U X)$$
where $P/U \cong \glr$ and $G/P \cong \PP^r$ is projective.
If $\overline{P \times_U X}$ is a $P$-equivariant projective
completion of $P \times_U X$ then $G \times_P (\overline{P \times_U X})$
is a projective completion of $G \times_U X$. When the action of
$U$ on $X$ extends to a $G$-action as above, we can choose
$\overline{P \times_U X}$ to be the closure of $P \times_U X$ in
$$\overline{G \times_U X} = \overline{G/U} \times X = \prrplus \times X;$$
that is,
$$\overline{P \times_U X} = \prr \times X \subseteq \overline{G \times_U X}$$
where $\prr = {\PP(\CC \oplus ((\CC^r)^* \otimes \CC^{r}))}$.
There is then a birational morphism
$$G \times_P (\overline{P \times_U X}) \to \overline{G \times_U X}$$
given by $[g,y] \mapsto gy$ which is an isomorphism over $G \times_U X$.
The resulting pullback $\hat{L} = \hat{L}^{(N)}$ to $G \times_P (\overline{P \times_U X})$
of $\calo_{\prrplus}(N) \otimes L$ is isomorphic to the induced
line bundle
$$G \times_P (\calo_{\prr} (N) \otimes L)$$
on $G \times_P (\overline{P \times_U X})$, where the $P$-action on 
$\calo_{\prr} (N) \otimes L$ is the restriction of the $G$-action on
$\calo_{\prrplus} (N) \otimes L$. 
If we regard  $G \times_P (\overline{P \times_U X})$ as a subvariety in
the obvious way of 
$$G \times_P (\overline{G \times_U X})= G \times_P (\prrplus \times X)$$
$$ \cong (G/P) \times \prrplus \times X \cong \PP^r \times \prrplus \times X$$
then the birational morphism 
$$G \times_P (\overline{P \times_U X}) \to \overline{G \times_U X} \cong \prr \times X$$
given by $[g,y] \mapsto gy$ extends to the projection
$$\PP^r \times \prrplus \times X \to \prrplus \times X$$
and so $\hat{L}^{(N)}$ is the restriction to $G \times_P (\overline{P \times_U X})$ of $\calo_{\prrplus}(N) \otimes L$.
Thus this line bundle $\hat{L}=\hat{L}^{(N)}$ is not
ample, but its tensor product $\hat{L}_\epsilon = \hat{L}_\epsilon^{(N)}$ with the pullback via the morphism
$$G \times_P (\overline{P \times_U X}) \to G/P \cong \PP^r,$$
of the
fractional line bundle $\calo_{\PP^r}(\epsilon)$,  where $\epsilon
\in \QQ \cap (0,\infty)$, provides an ample fractional linearisation
for the action of $G$ on $G \times_P ( \overline{P \times_U X})$ with, 
when $\epsilon$ is sufficiently small, an induced surjective birational
morphism
\begin{equation} \label{xhatu}
\widehat{\xu} =_{df} G \times_P ( \overline{P \times_U X}) /\!/_{\hat{L}_\epsilon} G
\to \overline{G \times_U X} /\!/ G = \xu
\end{equation}
(cf. \cite{K2,Reichstein}) which is an isomorphism over
$$(G \times_U X^{\bar{s}})/G \cong X^{\bar{s}}/U = X^s/U.$$
Note that $\hat{L}_\epsilon$ can be thought of as the bundle $G \times_P (\calo_{\prr} (N) \otimes L)$
on $G \times_P (\overline{P \times_U X})$, where now the $P$-action on 
$\calo_{\prr} (N) \otimes L$ is no longer the restriction of the $G$-action on
$\calo_{\prrplus} (N) \otimes L$ but has been twisted by $\epsilon$ times
the character of $P$ which restricts to the determinant on $\glr$.

\begin{rem}
It follows from  variation of GIT \cite{Ress} for the $G$-action
on $G \times_P (\overline{P \times_U X})$ that $\widehat{\xu}= \widehat{\xu}^{N,\e}$ 
is  independent of $N$ and $\epsilon$, provided that
$N$ is sufficiently large and $\epsilon >0$ is sufficiently small, depending on $N$.
\end{rem}

\begin{rem} When $\epsilon > 0$ the projective completion
$G \times_P (\overline{P \times_U X})$ equipped with the induced
ample fractional linearisation on $\hat{L}_\epsilon$ is not in general a 
reductive envelope for the $U$-action on $X$,
though it satisfies all the remaining conditions when (ii) is omitted from
Definition \ref{defn:separ} (cf. Remark \ref{rem:stablyfine}).  If we use the linearisation on $\hat{L}_0 
= \hat{L}$ instead, then we do obtain a reductive envelope, but it is not ample; nonetheless the
conditions of Proposition \ref{prop:main} are satisfied and 
we have
$$\overline{G \times_P(\overline{P \times_U X}) /\!/_{\hat{L}_0} G} = \overline{G \times_U X} /\!/G
= \xu.$$
\end{rem}

\begin{ex}
Let $U =\CC^+$ act linearly on a projective space $\PP^n$. Then we can choose
coordinates so that $1 \in {\rm Lie}(\CC^+) = \CC$ has Jordan normal 
form with blocks
$$ \left( \begin{array}{cccccccc}
0 & 1 & 0 & 0 & \cdots & 0\\
0 & 0 & 1 & 0 & \cdots & 0\\
&&& \cdots &&\\
0 & 0 & \cdots & 0 & 0 & 1\\
0 & 0 & \cdots & 0 & 0 & 0
\end{array} \right)$$
of sizes $k_1 + 1, \ldots, k_s + 1$ where $\sum_{j=1}^s (k_j +1) = n+1.$
Then the $\CC^+$ action extends to an action of $G= SL(2;\CC)$ via the
identifications 
$$ \CC^+ \cong \{ \left( \begin{array}{cccc}
1& a\\
0&1
\end{array} \right) : \ a \in \CC \} \leq G$$
and
$$\CC^{n+1} \cong \bigoplus_{j=1}^s {\rm Sym}^{k_j}(\CC^2)$$
where ${\rm Sym}^k(\CC^2)$ is the $k$th symmetric power of the
standard representation $\CC^2$ of $G = SL(2;\CC)$. Moreover
$$G/\CC^+ \cong \CC^2 \setminus \{ 0 \} \subseteq \CC^2 \subseteq
\PP^2 = \overline{G/\CC^+}$$ and thus we have
$$\PP^n /\!/ \CC^+ \cong {\rm Proj } (\CC[x_0, \ldots, x_n]^{\CC^+}) \cong (\PP^2 \times \PP^n)/\!/G$$
with respect to the linearisation $\calo_{\PP^2}(N) \otimes \calo_{\PP^n}(1)$
on $\PP^2 \times \PP^n$ for $N$ a sufficiently large positive integer.
When $G=SL(2;\CC)$ acts on $\PP^2$ we have
$(\PP^2)^{ss,G} = \CC^2$ (and $(\PP^2)^{s,G} = \emptyset$), so
since $N$ is large we have
$$(\PP^2 \times \PP^n)^{ss,G} \subseteq \CC^2 \times \PP^n
= (G \times_{\CC^+} \PP^n) \ \sqcup \  (\{ 0\} \times \PP^n)$$
and if semistability implies stability then
$$\PP^n /\!/ \CC^+ = (\PP^n)^{s,U}/\CC^+ \ \sqcup \ (\{0\} \times \PP^n)/\!/SL(2;\CC).$$
In this example the parabolic subgroup $P$ of $G = SL(2;\CC)$ is its Borel subgroup
$$ B = \{ \left( \begin{array}{cccc}
a & b\\
0& a^{-1}
\end{array} \right) : \  a \in \CC^*, b \in \CC \} $$
with $\overline{B/\CC^+} = \overline{\CC^*} = \PP^1$ and
$$\overline{B \times_{\CC^+} \PP^n} = \PP^1 \times \PP^n,$$
while $G \times_B \overline{B/\CC^+} = G \times_B \PP^1$ is the blow-up
of $\PP^2$ at the origin $0 \in \CC^2 \subseteq \PP^2$. Similarly
$G \times_B (\overline{B \times_{\CC^+} \PP^n})$ is the blow-up 
of $\overline{G \times_{\CC^+} \PP^n} \cong \PP^2 \times \PP^n$
along $\{ 0 \} \times \PP^n$, and its quotient $\widehat{\xu}$ is the blow-up of
$\PP^n /\!/\CC^+$ along its \lq boundary' 
$$\PP^n /\!/SL(2;\CC) \cong (\{0\} \times \PP^n)/\!/SL(2;\CC) \subseteq (\PP^2 \times \PP^n)/\!/SL(2;\CC)
= \PP^n/\!/\CC^+.$$
\end{ex}

Let us continue to assume that $U = \cplusr$ acts linearly on $X$ and that the action
extends to $G = \slrplus$. Notice that there are surjections
\begin{equation} \label{surj}
\overline{P \times_U X}^{ss,P,\epsilon} \to \widehat{\xu} \to \xu
\end{equation}
where $\overline{P \times_U X}^{ss,P,\epsilon}$ is the intersection of
$\overline{P \times_U X}$ with the $G$-semistable set in 
$G \times_P \overline{P \times_U X}$ with respect to the linearisation
$\hat{L}_\epsilon$, and $y_1,\ y_2 \in \overline{P \times_U X}^{ss,P,\epsilon}$
map to the same point in $\widehat{\xu}$ if and only if the closures of their
$P$-orbits $Py_1$ and $Py_2$ meet in $\overline{P \times_U X}^{ss,P,\epsilon}$.

Consider the linear action of the Levi subgroup $\glr \leq P$
on $\overline{P \times_U X} = \prr \times X$. It follows from the
Hilbert-Mumford criteria (Proposition \ref{sss} above) that 
\begin{equation} \label{pslss}
\overline{P \times_U X}^{ss,P,\epsilon} \subseteq 
\overline{P \times_U X}^{ss,\glr,\epsilon} \subseteq
\overline{P \times_U X}^{ss,SL(r;\CC)} \end{equation}
where $\overline{P \times_U X}^{ss,\glr,\epsilon}$
and $ \overline{P \times_U X}^{ss,SL(r;\CC)}$ (independent of
$\epsilon$) denote the $\glr$ and $SL(r;\CC)$-semistable sets
of $\overline{P \times_U X} $ after twisting the linearisation
by $\epsilon$ times the character $\det$ of $\glr$; this character
is of course trivial on $SL(r;\CC)$.

It is not hard to check that if the action of
$\glr$ on $\overline{P/U} = \prr$ is linearised with respect to
$\calo_{\prr}(1)$ by twisting by the fractional character $\frac{1}{2} \det$
then 
\begin{equation} \label{pmodu} 
\overline{P/U }^{ss,\glr,1/2} = \overline{P/U}^{s,\glr,1/2} = \glr \subseteq
(\CC^r)^* \otimes \CC^r \subseteq \prr.
\end{equation}
Thus, if instead of choosing $\epsilon$ close to 0 we choose $\epsilon$ to
be approximately $N/2$, where $N$ is the sufficiently large positive
integer chosen above, then we see from the Hilbert-Mumford criteria
(Proposition \ref{sss}) that
$$\overline{P \times_U X}^{ss,\glr,\epsilon} =
(\prr \times X)^{ss,\glr,\epsilon} = \glr \times X$$
and so quotienting we get
$$\overline{P \times_U X}/\!/_{\hat{L}^{(N)}_{N/2}} \glr = X.$$
A GIT quotient of a nonsingular complex projective variety $Y$ by a linear action of $\glr$ can always be constructed
by first quotienting by $SL(r;\CC)$ and then quotienting by the
induced linear action of $\CC^* = \glr/SL(r;\CC)$: we have
$$Y/\!/ \glr = (Y/\!/ SL(r;\CC)  ) \ /\!/ \ \CC^*.$$
 Therefore if we set
\begin{equation} \label{definecalx}
\calx = \overline{P \times_U X}/\!/_{\hat{L}^{(N)}} SL(r;\CC)
= (\prr \times X)/\!/_{\hat{L}^{(N)}} SL(r;\CC)
\end{equation}
for $N>0$ sufficiently large, then $\calx$ is a projective variety with
a linear action of $\CC^*$ which we can twist by $\epsilon$ times the
standard character of $\CC^*$, such that when $\epsilon = N/2$ we get
\begin{equation} \label{calxquot}
\calx /\!/_{N/2} \CC^* \ \cong \  X
\end{equation}
while for $\epsilon >0$ sufficiently small it follows from
(\ref{surj}) and (\ref{pslss}) that we have a surjection
from an open subset of $\calx /\!/_{\epsilon} \CC^*$ onto
$\widehat{\xu}$, and hence onto $\xu$. More precisely, the 
inclusion
$$(\hat{\calo}_{\hat{L}_{\epsilon}}(G \times_P (\overline{P \times_U X})))^G
= (\hat{\calo}_{\hat{L}_{\epsilon}}(\overline{P \times_U X}))^P
\subseteq (\hat{\calo}_{\hat{L}_{\epsilon}}(\overline{P \times_U X}))^{\glr}$$
induces a rational map 
\begin{equation} \label{ratmap}
\calx/\!/_{\epsilon} \CC^* = \overline{P \times_U X}/\!/_{\hat{L}_{\epsilon}} \glr
-\ -\ \to G \times_P \overline{P \times_U X}/\!/_{\hat{L}_{\epsilon}} G =
\widehat{\xu} \end{equation}
whose composition with the surjection
$$\overline{P \times_U X}^{ss,\glr,\epsilon} \to \calx /\!/_{\epsilon} \CC^*$$
induced by the inclusion
$$(\hat{\calo}_{\hat{L}_{\epsilon}}(\overline{P \times_U X}))^{\glr}
\subseteq (\hat{\calo}_{\hat{L}_{\epsilon}}(\overline{P \times_U X}))^{SL(r;\CC)}$$
is the rational map
$$\overline{P \times_U X}^{ss,\glr,\epsilon} - \ - \ \to 
G \times_P \overline{P \times_U X}/\!/_{\hat{L}_{\epsilon}} G =
\widehat{\xu}$$
which restricts to a surjection
$$\overline{P \times_U X}^{ss,P,\epsilon} \to
\widehat{\xu}.$$
Hence the restriction of $\calx/\!/_{\epsilon} \CC^* 
-\ -\ \to \widehat{\xu}$ to its domain of definition is surjective.

\begin{defn} \label{defn4p4}
Let $(\calx/\!/_{\epsilon} \CC^*)^{\hat{ss}}$ denote the open
subset of $\calx/\!/_{\epsilon} \CC^*$ which is the domain of
definition of the rational map (\ref{ratmap}) from
$\calx/\!/_{\epsilon} \CC^*$ to $\widehat{\xu}$, where as above
$\calx$ is the projective variety
$$\calx = (\overline{P \times_U X})/\!/ SL(r;\CC)$$
with the induced linear $\CC^*$-action, and $0<\epsilon <\!< 1$. 
Let $(\calx/\!/_{\epsilon} \CC^*)^{\hat{s}}$ be the open
subset 
$\overline{P \times_U X}^{s,P,\epsilon}/\glr$ of
$$\overline{P \times_U X}^{s,\glr,\epsilon}/\glr = 
(\overline{P \times_U X}^{s,\glr,\epsilon}/SL(r;\CC))/\CC^*$$
$$
= \calx^{s,\epsilon}/\CC^* \subseteq
\calx/\!/_{\epsilon} \CC^*.$$
Let $\calx^{\hat{ss},\epsilon} = \pi^{-1}((\calx/\!/_{\epsilon} \CC^*)^{\hat{ss}})$
and $\calx^{\hat{s},\epsilon} = \pi^{-1}((\calx/\!/_{\epsilon} \CC^*)^{\hat{s}})$
where $\pi: \calx^{ss,\epsilon} \to \calx/\!/_{\epsilon} \CC^*$ is the quotient
map, so that
$$(\calx/\!/_{\epsilon} \CC^*)^{\hat{s}} = \calx^{\hat{s},\epsilon}/\CC^*.$$
\end{defn}

Then we have

\begin{prop} \label{xquote}
If $\epsilon>0$ is sufficiently small, the rational map from
$\calx/\!/_{\epsilon} \CC^*$ to $\widehat{\xu}$ induced by the inclusion
of $ 
 (\hat{\calo}_{\hat{L}_{\epsilon}}(\overline{P \times_U X}))^P$ in
$ (\hat{\calo}_{\hat{L}_{\epsilon}}(\overline{P \times_U X}))^{\glr}$
restricts to surjective morphisms
$$(\calx/\!/_{\epsilon} \CC^*)^{\hat{ss}} \to \widehat{\xu} \to \xu$$
and
$$(\calx/\!/_{\epsilon} \CC^*)^{\hat{s}} \to X^s/U.$$
\end{prop}

\begin{rem} \label{rem4p6}
Using the theory of variation of GIT \cite{DolgHu, Ress, Thad}, as described in Remark \ref{flips},
we can relate the quotient $\calx/\!/_{\epsilon} \CC^*$ which appears in Proposition
\ref{xquote} above to $\calx/\!/_{N/2} \CC^* \cong X$ via a sequence of flips which occur as walls
are crossed between the linearisations corresponding to $\epsilon$ and to $N/2$. Thus we have a
diagram
$$ \begin{array}{ccccccccc}
(\calx/\!/_{\epsilon} \CC^*)^{\hat{s}} & \subseteq & (\calx/\!/_{\epsilon} \CC^*)^{\hat{ss}} &
\subseteq & \calx/\!/_{\epsilon} \CC^* & \leftarrow - \rightarrow &
X = \calx/\!/_{N/2} \CC^*  \\
\downarrow & &  \downarrow & & & \mbox{flips} & \\
X^{{s}}/U & \subseteq &
\widehat{\xu} & & & & \\
|| & & \downarrow & & & & \\
X^s/U & \subseteq & \xu & & & & 
\end{array}
$$
where the vertical maps are all {\em surjective}, in contrast to (\ref{xudiagram}), and
the inclusions are all open.

Note also that by variation of GIT if $0<\epsilon <\!< 1$ there is a birational
surjective morphism
$$\calx/\!/_{\epsilon} \CC^* \to \calx/\!/_{0} \CC^*.$$
When $\epsilon = 0$  the 
inclusion
$$(\hat{\calo}_{\hat{L}_{0}}(G \times_P (\overline{P \times_U X})))^G
= (\hat{\calo}_{\hat{L}_{0}}(\overline{P \times_U X}))^P
\subseteq (\hat{\calo}_{\hat{L}_{0}}(\overline{P \times_U X}))^{\glr}$$
induces a rational map 
\begin{equation} \label{ratomap}
\calx/\!/_{0} \CC^* = \overline{P \times_U X}/\!/_{\hat{L}_{0}} \glr
-\ -\ \to G \times_P \overline{P \times_U X}/\!/_{\hat{L}_{0}} G =
{\xu} \end{equation}
whose composition with the surjective morphism
$$\calx/\!/_{\epsilon} \CC^* \to \calx/\!/_{0} \CC^*$$
is the composition of (\ref{ratmap}) with the surjective morphism 
$\widehat{\xu} \to \xu$. Thus the restriction of the rational map 
(\ref{ratomap}) from $\calx/\!/_{0} \CC^*$ to ${\xu}$ to its domain
of definition is surjective.
\end{rem}

\begin{rem}
\label{remwonderful}
Note that the GIT quotient $\prr /\!/SL(r;\CC)$ is isomorphic to $\PP^1$, but we
do not have $(\prr)^{ss} = (\prr)^s$ for this action of $SL(r;\CC)$. It
is therefore convenient to replace the compactification $\prr$ of $\glr$ by its
{\em wonderful compactification} $\widetilde{\prr}$ given by blowing up
$\prr= \{[z:(z_{ij})_{i,j=1}^r]\}$ along the (proper transforms of the) subvarieties 
defined by
$$z=0 \mbox{ and } {\rm rank}(z_{ij}) \leq \ell$$
for $\ell = 0,1,\ldots,r$ and by
$${\rm rank}(z_{ij}) \leq \ell$$
for $\ell = 0,1,\ldots,r-1$ \cite{Kausz}. The action of $SL(r;\CC)$ on
$\widetilde{\prr}$, linearised with respect to a small perturbation $\calo_{\widetilde{\prr}}(1)$
of the pullback of $\calo_{\prr}(1)$, satisfies
$$\widetilde{\prr}^{ss} = \widetilde{\prr}^s \mbox{ and } \widetilde{\prr}/\!/SL(r;\CC) \cong \PP^1.$$
Thus if we replace $\overline{P \times_U X}=\prr \times X$ with 
$$\widetilde{{P \times_U X}} = \widetilde{\prr} \times X$$
 and define
$\widetilde{\xu} = G \times_P(\widetilde{P \times_U X})/\!/_{\hat{L}_\epsilon} G$ and 
\begin{equation} \label{newdefncalx}
\widetilde{\calx} =  \widetilde{P \times_U X}/\!/_{\hat{L}^{(N)}} SL(r;\CC)
= (\widetilde{\prr} \times X)/\!/_{\hat{L}^{(N)}} SL(r;\CC)
\end{equation}
for $N >\!> 0$, then all the properties of $\calx$ given above still hold for $\widetilde{\calx}$,
and in addition $\widetilde{\calx}$ fibres over $\PP^1$ as
$$\widetilde{\calx} 
= (\widetilde{\prr} \times X)/\!/_{\hat{L}^{(N)}} SL(r;\CC)
= (\widetilde{\prr}^{ss} \times X)/ SL(r;\CC) $$
$$ \to \  \  \  
\widetilde{\prr}^{ss}/ SL(r;\CC) = \PP^1$$
with fibres isomorphic to the quotient of $X$ by the finite centre of $SL(r;\CC)$.
If $X$ is nonsingular then it turns out that $\widetilde{\calx}$ and $\widetilde{\calx} /\!/_{\epsilon} \CC^*$ (for $0 < \epsilon
<\!< 1$) and
$\widetilde{\xu}$ are orbifolds, so that $\widetilde{\xu}$ is a projective completion
of $X^s/U$ which is a partial desingularisation 
of $\xu$ (cf. Remark 4.1).
\end{rem}

\subsection{General $\cplusr$ actions}
Of course the constructions described in $\S$4.1 only work if the action
of $U = \cplusr$ on $X$ extends to an action of $G = SL( \CC \oplus \lieu)$, which is a
rather special situation when the ring of invariants $\hat{\calo}_L(X)^U$ is 
always finitely generated. Moreover at least {\em a priori} these
constructions may depend on the choice of
this extension, although
$\overline{G \times_U X}/\!/G = \xu = {\rm Proj}(\hat{\calo}_L(X)^U)$ depends
only on the linearisation of the $U$-action on $X$. 
So next we need to consider what happens if the linear $U$-action on $X$ does not extend to
a linear action of $G$.
Suppose that we
can associate to the linear $U$-action on $X$ a normal projective variety $Y$  containing $X$,
with an action of $G = SL( \CC \oplus \lieu)$
and a $G$-linearisation on a line bundle $L_Y$, which restricts to the given 
linearisation of the $U$-action on $X$ and is such that every $U$-invariant in a
finite fully separating set of $U$-invariants on $X$ extends to a $U$-invariant on $Y$.
Then we can embed $X$ in the $G$-variety
$$\prrplus \times Y$$ as $\{\iota \} \times X$ where
$\iota \in (\CC^r)^* \otimes \CC^{r+1} \subseteq \prrplus$
is the standard embedding of $\CC^r$ in $\CC^{r+1}$, and the
closure of $GX \cong G \times_U X$ in $\prrplus \times Y$ will provide
us with a reductive envelope $\overline{G \times_U X}$.
Therefore
we will next consider how, given any linearised $U$-action on $X$, we can
choose a $G$-variety $Y$ with these properties. We will find that for any sufficiently divisible
positive integer $m$ we can choose such a variety $Y_{m}$ in a canonical
way, depending only on $m$ and the linear action of $U$ on $X$, giving us
a reductive envelope $\overline{G \times_U X}^{m}$.

Let $S$ be any finite fully separating set of invariants (in the sense of Definition
\ref{defn:separ}) on $X$. By replacing the elements of $S$ with suitable
powers of themselves, we can assume that 
$S \subseteq H^0(X,L^{\otimes m})^U$
for some $m>0$ for which
$L^{\otimes m}$ is very ample. Then $X \subseteq \PP(H^0(X,L^{\otimes m})^*)$ and
every $\sigma \in S$ extends to a $U$-invariant section of $\calo(1)$ on
$\PP(H^0(X, L^{\otimes m})^*)$.

Now consider the linear action of $U$ on $V_m = H^0(X,L^{\otimes m})^*$,
and let $P$ be the parabolic subgroup of $G = \slrplus$ with unipotent 
radical $U$, as at (\ref{parabolic}) above. Since $P$ is
a semi-direct product
$$P = U \rtimes \glr$$
we have
$$P \times_U V_m \cong \glr \times V_m$$
with the $P$-action on $\glr \times V_m$ given for $(h,v) \in \glr \times V_m$ by
$$p.(h,v) = (gh, (h^{-1}uh).v)$$
where $p=gu$ with $g \in \glr$ and $u \in U$, and $h^{-1}uh$ acts on
$v \in V_m$ via the given $U$-action. Of course $\glr \times V_m$ is an affine variety
with
$$\calo(\glr \times V_m) \cong \CC[h_{ij}, (\det h)^{-1}, v_k]$$
where $\det h$ is the determinant of the $r \times r$ matrix $(h_{ij})_{i,j=1}^r$ and
$(v_k)$ are coordinates on $V_m$. Let
\begin{equation} \label{phivm}
\phi_{V_m}: \CC^r = {\rm Lie}\ U \to {\rm Lie}(GL(V_m))
\end{equation}
be the infinitesimal action of $U$ on $V_m$ and let $U_{V_m}$ be its image in 
${\rm Lie}(GL(V_m))$. Since $U$ is unipotent we have
$$V_m \supseteq U_{V_m}(V_m) \supseteq (U_{V_m})^2(V_m) \supseteq \cdots \supseteq 
(U_{V_m})^{\dim V_m}(V_m) =  0 $$
where $$
(U_{V_m})^j(V_m)   =  \{u_1 u_2 \cdots u_j(v) \ : \ u_1,\ldots,u_j \in U_{V_m}, \ v \in V_m \} $$
$$  =  \{\phi_{V_m}(\tilde{u}_1) \phi_{V_m}(\tilde{u}_2) \cdots 
\phi_{V_m}(\tilde{u}_j)(v) \ : \ \tilde{u}_1,\ldots,\tilde{u}_j \in {\rm Lie}\ U, \ v \in V_m \}.
$$
For $0 \leq j \leq \dim V_m - 1$ let $\Theta_{j,m}$ be the complex vector space consisting of
all polynomial functions 
$$\theta: (\CC^r)^j \times ((\CC^r)^* \otimes \CC^r)^{\dim V_m -1} \to \CC$$
$$ (u_1,\ldots,u_j,h_1,\ldots,h_{\dim V_m -1}) \mapsto \theta(u_1,\ldots,u_j,h_1,\ldots,h_{\dim V_m -1})$$
which are simultaneously homogeneous of degree 1 in the coordinates of each $u_i \in \CC^r$ separately,
 for $1 \leq i \leq j$,
and homogeneous of total degree $r(\dim V_m -1) - (r-1)j$ in the coordinates of all the $h_k \in (\CC^r)^* \otimes \CC^r$ together, for $1 \leq k \leq \dim V_m -1$. Let
\begin{equation} \label{wm}
W_m =  \bigoplus_{j=0}^{\dim V_m -1} 
\Theta_{j,m} \otimes (U_{V_m})^j(V_m).
\end{equation}
Then we can embed $V_m$ linearly into 
$W_m$    via
\begin{equation} \label{defnpsi} v \mapsto \psi(v) =  \sum_{j=0}^{\dim V_m -1} \psi_j(v) \end{equation}
where  $\psi_j(v):(\CC^r)^{j}\times((\CC^r)^*  
\otimes\CC^r)^{(\dim V_m - 1)}\to (U_{V_m})^j(V_m)$
for $0 \leq j \leq \dim V_m -1$ and $v \in V_m$ sends
$ (u_1 , \ldots , u_j, h_1 ,\ldots , h_{\dim V_m -1})$ to the product of 
$$  
\det(h_1) \det (h_2) \cdots \det (h_{\dim V_m - j - 1})$$
with $$   \phi_{V_m}(h_{\dim V_m -1}u_1) \cdots
\phi_{V_m}(h_{\dim V_m -j }u_j)(v) .$$
In particular 
$\psi_0(v): ((\CC^r)^*  
\otimes\CC^r)^{(\dim V_m - 1)}) \to (U_{V_m})^0 (V_m) = V_m$ sends
$(h_1 ,\ldots , h_{\dim V_m -1})$
to 
$$( \det(h_1) \det (h_2) \cdots \det (h_{\dim V_m  - 1})   ) v 
\ \in \ V_m
. $$ This embedding of $V_m$ into $W_m$ is $U$-equivariant
with respect to the linear $U$-action on $W_m$
for which the infinitesimal action of $u \in {\rm Lie} \ U = \CC^r$ is given by
\begin{equation} \label{uaction}
u \ ( \sum_{j=0}^{\dim V_m -1} \alpha_j) =  \sum_{j=0}^{\dim V_m - 1} u \cdot \alpha_{j+1}\end{equation}
with $u \cdot \alpha_{j+1}:(\CC^r)^j \times ((\CC^r)^* \otimes \CC^r)^{\dim V_m -1} \to (U_{V_m})^j(V_m)
 $ defined to be 0 for
$j=\dim V_m$ and for any
$0 \leq j < \dim V_m$ 
defined by 
$$
u\cdot \alpha_{j+1}(u_1, \ldots , u_j, h_1 ,\ldots , h_{\dim V_m -1 })$$
$$ = \alpha_{j+1}( u_1 
, \ldots , u_j , {\rm adj}(h_{\dim V_m - 1 - j})(u),h_1, \ldots,h_{\dim V_m -1}). $$
Here ${\rm adj}(h)$ is the adjoint matrix of $h$ (so that $h \ {\rm adj}(h) = {\rm adj}(h) \ h$
is $\det (h)$ times the identity matrix $\iota$), which is homogeneous of degree $r-1$ in
the coefficients of $h$. 
Moreover this linear action of $U$ on $W_m$ extends to a
linear action of $P=U \rtimes \glr$ on $ 
W_m$ where
$g \in \glr$ acts as
\begin{equation} \label{glraction}
g \ ( \sum_{j=0}^{\dim V_m -1} \alpha_j) =  \sum_{j=0}^{\dim V_m -1}(\det g)^{j}\ g \cdot \alpha_j
\end{equation}
with $g\cdot \alpha_j:(\CC^r)^j \times ((\CC^r)^* \otimes \CC^r)^{\dim V_m -1} \to (U_{V_m})^j(V_m)
 $ defined for any
$0 \leq j < \dim V_m$ and any $ \alpha_j:(\CC^r)^j \times ((\CC^r)^* \otimes \CC^r)^{\dim V_m -1} \to (U_{V_m})^j(V_m)$ by 
$$g\cdot \alpha_j(   
u_1 , \ldots , u_j, h_1 ,\ldots , h_{\dim V_m -1 }) $$
$$ =  \alpha_j(g u_1 
, \ldots , g u_j, g h_1 g^{-1}, \ldots, g h_{\dim V_m - 1}g^{-1}).$$

Since the $U$-action on $\PP(W_m)$  
extends to a linear $P$-action, we can
construct the projective variety 
$$Y_m = G \times_P \PP(W_m) $$ 
and we can equip $Y_m$ with the line bundles $L_{Y_m}^{\epsilon} = G \times_P \calo_{\PP(W_m)}(1)$ for $\epsilon \in \QQ$ where
the $P$-action on $\calo_{\PP(W_m)}(1)$ is induced from the given $P$-action on $W_m$ twisted
by $\epsilon$ times the pullback to $P$ of the character $\det$ of $P/U \cong \glr$. Equivalently
$L_{Y_m}^{\epsilon}$ is the tensor product of $L^0_{Y_m}$ with the pullback via
$$G \times_P \PP(W_m) \to G/P \cong \PP^r$$
of $\calo_{\PP^r}(\epsilon)$. 

\begin{rem} \label{pre4.9}
The action of $P$ on $W_m$ extends to a linear action of $G=\slrplus$ on the vector space
$\calw_m \supseteq W_m$, where
$$\calw_m = 
H^0(G \times_P \PP(W_m \otimes (\bigoplus_{\lambda \in \Delta} \CC_\lambda));G \times_P \calo(1))$$
for a sufficiently large finite set $\Delta$ of weights of $P$, with $\CC_\lambda$ denoting a 
copy of $\CC$ on which $P$ acts with weight $\lambda$.

Thus $ G \times_P \PP(W_m) \subseteq G \times_P \PP(\calw_m) \cong G/P \times \PP(\calw_m) \cong \PP^r \times \PP(\calw_m)$, and with respect to this
identification we have
$$L^0_{Y_m} = \calo_{\PP(\calw_m)}(1)|_{Y_m}.$$
Hence if $\e >0$ it follows that $L^\e_{Y_m}$ is the restriction  of
the line bundle 
$\calo_{\PP^r}(\e) \otimes \calo_{\PP(\calw_m)}(1)$ on $\PP^r \times 
\PP(\calw_m)$, and so $L^\e_{Y_m}$ is ample.
\end{rem}

This gives us for every sufficiently divisible positive integer $m$ a $U$-equivariant embedding
$$X \subseteq \PP(H^0(X,L^{\otimes m})^*) = \PP(V_m) \subseteq \PP(W_m) $$
in a projective variety with a linear $P$-action,
such that every $\sigma$ in the finite fully separating set of invariants 
$S$ extends to a $U$-invariant linear functional
on $V_m$, and hence extends to a $U$-invariant (in fact $P$-invariant) linear functional on $W_m$ defined
by $\sum_{j=0}^{\dim V_m} \alpha_j \mapsto \sigma(\alpha_0(\iota,\ldots,\iota))$.
As each $\sigma \in S$ extends to a $P$-invariant linear functional
on $W_m$, it extends to a $G$-invariant section of $L^0_{Y_m}$. Thus from $\S$4.1 we have
the following proposition.

\begin{prop} \label{mredenv}
Let $X$ be embedded in $Y_m = G \times_P \PP(W_m)$ as above, for a sufficiently divisible
positive integer $m$, and let $\iota \in (\CC^r)^* \otimes \CC^{r+1} \subseteq \prrplus$
be the standard embedding of $\CC^r$ in $\CC^{r+1}$. If $N$ is sufficiently large (depending on $m$), then the linear action of $U = \cplusr$
on $X$ has a reductive envelope given by the closure 
$\overline{G \times_U X}^m$ of $G \times_U X$ embedded in $\prrplus \times Y_m$ as 
$G(\{\iota\} \times X)$, equipped with the restriction of the $G$-linearisation on $\calo_{\prrplus}(N) \otimes L^0_{Y_m}$.
\end{prop}


Note that by Remark \ref{pre4.9} the line bundle $L^0_{Y_m}$ is not in general ample (although $L^{\epsilon}_{Y_m}$ is ample for any $\epsilon > 0$),
 so we do not necessarily have an ample reductive envelope here. 
Nonetheless by Proposition \ref{prop:main} if $X^{\overline{s}}$ and $X^{\overline{nss}}$ are defined using this
reductive envelope we have

\begin{corollary} \label{march}
$X^{\overline{s}} \subseteq X^s \subseteq X^{ss} \subseteq X^{\overline{nss}}.$
\end{corollary}

If moreover the ring of invariants $\hat{\calo}_L(X)^U$ is finitely generated and $m$ is sufficiently
divisible that $\hat{\calo}_{L^{\otimes m}}(X)^U$ is generated by $H^0(X,L^{\otimes m})^U$,
then for $N>\!>0$ the restriction map
$$\rho_m: \bigoplus_{k \geq 0} H^0(\overline{G \times_U X}^m,(\calo_{\prrplus}(N) \otimes L^0_{Y_m})^{\otimes k})^G
\to \hat{\calo}_{L^{\otimes m}}(X)^U$$
is an isomorphism and
$\xu = \Proj \hat{\calo}_L(X)^U$ is the canonical projective completion of $\overline{G \times_U X}^m/\!/G$
(see Proposition \ref{prop:main} again).
Even when the ring of invariants $\hat{\calo}_L(X)^U$ is not finitely generated, if $m$ is sufficiently
divisible that $H^0(X,L^{\otimes m})$ contains a finite fully separating set of invariants, then
for any multiple $m'=k'm$ of $m$
the subalgebra $$\hat{\calo}_L^{m'}(X)^U$$ of $\hat{\calo}_L(X)^U$ generated by $H^0(X,L^{\otimes m'})^U$
is finitely generated and provides a projective completion 
$$\overline{\xu}^m = \Proj \hat{\calo}_L^{m'}(X)^U$$
of $\xu$, while
the restriction of $\rho_m$ to the subalgebra  of $$\bigoplus_{k \geq 0} H^0(\overline{G \times_U X}^m,
(\calo_{\prrplus}(N) \otimes L^0_{Y_m})^{\otimes k})^G$$ generated by
$H^0(\overline{G \times_U X}^m,(\calo_{\prrplus}(N) \otimes L^0_{Y_m})^{\otimes k'})^G$
gives an isomorphism onto $\hat{\calo}_L^{m'}(X)^U$.

By analogy with the construction of $\widehat{\xu}$ and $\widetilde{\xu}$ in $\S$4.1, let us also 
consider the closures $\widehat{G \times_U X}^m$ and $\widetilde{G \times_U X}^m$
of $G \times_U X = G(\{\iota\} \times X)$ in $G \times_P (\prr \times Y_m)$ and
$G \times_P (\widetilde{\prr} \times Y_m)$ respectively. Since $Y_m = G \times_P \PP(W_m)$ we have
$$\widehat{G \times_U X}^m \cong G \times_P (\overline{P \times_U X}^m) \cong \overline{G \times_U X}^m$$
 and 
$$\widetilde{G \times_U X}^m \cong G \times_P (\widetilde{P \times_U X}^m)$$
where $\overline{P \times_U X}^m$ and $\widetilde{P \times_U X}^m$ are the closures
of $P \times_U X= P(\{\iota\} \times X)$ in $\prr \times \PP(W_m)$ and $\widetilde{
\prr} \times \PP(W_m)$ respectively.

\begin{defn}
Let $\hat{L}_\epsilon = \hat{L}_\epsilon^{(N)}$ be the tensor product of the pullback
via
$$\widehat{G \times_U X}^m \cong G \times_P (\overline{P \times_U X}^m) \to G/P \cong \PP^r$$
of $\calo_{\PP^r}(\epsilon)$ with the line bundle
$G \times_P \hat{L}^{(N)}$ 
on $ G \times_P (\overline{P \times_U X}^m)$
where $$\hat{L}^{(N)} =
 \calo_{{\prr}}(N) \otimes \calo_{\PP(W_m)}(1)|_{\overline{P \times_U X}^m};$$
 equivalently $\hat{L}_\epsilon^{(N)} = G \times_P \hat{L}^{(N)}$ where
 the action of $P$ on $\hat{L}^{(N)}$ is twisted by $\e$ times the character
 $\det$ of $P/U \cong \glr$. Similarly
let $\tilde{L}_\epsilon = \tilde{L}_\epsilon^{(N)}$ be the tensor product of the pullback
via
$$\widetilde{G \times_U X}^m \cong G \times_P (\widetilde{P \times_U X}^m) \to G/P \cong \PP^r$$
of $\calo_{\PP^r}(\epsilon)$ with the line bundle
$G \times_P \tilde{L}^{(N)}$ 
on $ G \times_P (\widetilde{P \times_U X}^m)$, where $\tilde{L}^{(N)} =
 \calo_{\widetilde{\prr}}(N) \otimes \calo_{\PP(W_m)}(1)|_{\widetilde{P \times_U X}^m}$ and as in Remark \ref{remwonderful} 
$\calo_{\widetilde{\prr}}(1)$ is a small perturbation of the pullback of $\calo_{{\prr}}(1)$
along $\widetilde{\prr} \to \prr$ such that the $SL(r;\CC)$-action lifts to 
$\calo_{\widetilde{\prr}}(1)$ and satisfies $(\widetilde{\prr})^{ss} = (\widetilde{\prr})^s$
and $\widetilde{\prr}/\!/SL(r;\CC) \cong \PP^1$.
Note that the line bundles $\hat{L}_\epsilon^{(N)}$ on
$\widehat{G \times_U X}^m \cong G \times_P (\overline{P \times_U X}^m)$
and $\calo_{\prrplus}(N) \otimes L^{\epsilon}_{Y_m}$ on 
$$\widehat{G \times_U X}^m \cong G \times_P (\overline{P \times_U X}^m) \cong \overline{G \times_U X}^m$$
are both $G$-invariant and both restrict to $\calo_{{\prr}}(N) \otimes \calo_{\PP(W_m)}(1)$
on $\overline{P \times_U X}^m$ with the same $P$-action, so they are isomorphic to each other.
\end{defn}

The line bundles $\hat{L}_\epsilon = \hat{L}_\epsilon^{(N)}$ and
 $\tilde{L}_{\epsilon} = \tilde{L}_{\epsilon}^{(N)}$  on $\overline{G \times_U X}^m$ 
 and $\widetilde{G \times_U X}^m$ are ample for 
$\epsilon >0$, and the $G$-actions on $\overline{G \times_U X}^m$ 
 and 
$\widetilde{G \times_U X}^m$ lift to  linear actions on $\hat{L}_{\epsilon}$ and
$\tilde{L}_{\epsilon}$.

\begin{defn} \label{calxm}
For a positive integer $N$ sufficiently large
(depending on $m$) 
let
$$ {\calx}_m = \overline{P \times_U X}^m /\!/_{\hat{L}^{(N)}} SL(r;\CC)
\ \ \mbox{ and } \ \ \tilde{\calx}_m = \widetilde{P \times_U X}^m /\!/_{\tilde{L}^{(N)}} SL(r;\CC)$$
and for $\epsilon >0$ sufficiently small (depending on $m$ and $N$) let
$$ \widehat{\xu}^{m}  
= \overline{G \times_U X}^m /\!/_{\hat{L}_\epsilon^{(N)}} G
\ \ \mbox{ and } \ \ \widetilde{\xu}^{m} 
 = \widetilde{G \times_U X}^m /\!/_{\tilde{L}_\epsilon^{(N)}} G.$$
\end{defn}

\begin{rem} The line bundles $\hat{L}^{(N)}$ and $\tilde{L}^{(N)}$ are
ample on $\overline{P \times_U X}^m$ and $\widetilde{P \times_U X}^m$,
and for $\e >0$ the line bundles $\hat{L}^{(N)}_\e $ and $\tilde{L}^{(N)}_\e $ are
ample on $\overline{G \times_U X}^m$ and $\widetilde{G \times_U X}^m$. Thus 
it follows from  variation of GIT \cite{Ress} for the $SL(r;\CC)$-actions
on $\overline{P \times_U X}^m$ and $\widetilde{P \times_U X}^m$
and the $G=SL(r+1;\CC)$-actions
on $\overline{G \times_U X}^m$ and $\widetilde{G \times_U X}^m$ that 
${\calx}_m$ and 
$ \tilde{\calx}_m$ and $\widehat{\xu}^m$ and $\widetilde{\xu}^m$ are
 independent of $N$ and $\e $, provided that
$m$ is fixed and $N$ is sufficiently large,
depending on $m$, and $\e$ is sufficiently small, depending on $N$ and $\e$.
\end{rem}

\begin{rem} \label{bbb}
Recall that $\hat{L}_\epsilon^{(N)}$ can be identified with the line bundle
$G \times_P \hat{L}^{(N)}$ 
on $ G \times_P (\overline{P \times_U X}^m)$ when the action of $P$ on
$\hat{L}^{(N)}$ is twisted by $\e$ times the character
 $\det$ of $P/U \cong \glr$. 
The character $\det$ extends to a section of the line bundle
$\calo_{\prr}(r)$ over the projective completion $\prr$ of $\glr$
on which $P$ acts by multiplication by the character $\det$, and thus to
a section of the line bundle $G \times_P \calo_{\prr}$ over $ G \times_P (\overline{P \times_U X}^m)$ when the action of $P$ on $\calo_{\prr}$ is
twisted by this character.
Tensoring with this section gives us an injection
 \begin{equation} \label{klob} H^0(\overline{G \times_U X}^m,(\hat{L}_0^{(N)})^{\otimes k})^G
 \to H^0(\overline{G \times_U X}^m,(\hat{L}_\e ^{(N + r\e )})^{\otimes k})^G \end{equation}
 where $\e = 1/k$, whose composition with the injection given by the restriction
map 
$$\rho_m^\e: \bigoplus_{k \geq 0} H^0(\overline{G \times_U X}^m,(\hat{L}_\e^{(N + r\e)})^{\otimes k})^G
\to \hat{\calo}_{L^{\otimes m}}(X)^U$$
is the restriction map $\rho_m=\rho_m^0$. If 
 the ring of invariants $\hat{\calo}_L(X)^U$ is finitely generated and
 $m$ is sufficiently
divisible that $\rho_m$ is an isomorphism, 
it follows
 that (\ref{klob}) is an isomorphism and thus that
$$(\overline{G \times_U X}^m)^{ss,\e,G} \subseteq (\overline{G \times_U X}^m)^{ss,0,G} $$
and that the inclusion of $(\overline{G \times_U X}^m)^{ss,\e,G}$
in $(\overline{G \times_U X}^m)^{ss,0,G}$ induces a birational surjective morphism
$$\widehat{\xu}^m \to \xu = \Proj(\hat{\calo}_L(X)^U).   
$$
Even when $\hat{\calo}_L(X)^U$ is not finitely generated,
if $m$ and $k$ are sufficiently divisible that $H^0(X,L^{\otimes m})^U$
contains a fully separating set of invariants and 
$\calo_{(\hat{L}_\e^{(N + r\e)})^{\otimes k}}(\overline{G \times_U X}^m)$
is generated by
$H^0(\overline{G \times_U X}^m,(\hat{L}_\e^{(N + r\e)})^{\otimes k})^G$, 
then we
have a birational surjective morphism
$$\widehat{\xu}^m \to \Proj(\hat{\calo}_L^{m'}(X)^U) = \overline{\xu}^{m'}$$
where $m'=km$. The same is also true when $\widehat{\xu}^m$ is replaced with
$\widetilde{\xu}^m$.
\end{rem}

As in Proposition \ref{xquote} and Remark \ref{remwonderful} we obtain

\begin{theorem} \label{genaction}
If $m$ is a sufficiently divisible positive integer and $N>\!>0$ then $\calx_{m}$ and
$\widetilde{\calx}_m$ are projective varieties
with linear actions of $\CC^*$ which we can twist by $\epsilon$ times the standard character of
$\CC^*$, such that when $\epsilon = N/2$ we have
$$\calx_{m}/\!/_{N/2} \CC^* \cong X \ \mbox{ and } \ 
\widetilde{\calx}_{m}/\!/_{N/2} \CC^* \cong X,$$
while if $\epsilon >0$ is sufficiently small then the rational maps from
$\calx_{m}/\!/_{\epsilon} \CC^*$ to $\widehat{\xu}^{m}$ and from
$\widetilde{\calx}_{m}/\!/_{\epsilon} \CC^*$ to $\widetilde{\xu}^{m}$
induced by the inclusions
of $ 
 (\hat{\calo}_{\hat{L}_{\epsilon}}(\overline{P \times_U X}^{m}))^P$ in
$ (\hat{\calo}_{\hat{L}_{\epsilon}}(\overline{P \times_U X}^{m}))^{\glr}$ and
$  (\hat{\calo}_{\hat{L}_{\epsilon}}(\widetilde{P \times_U X}^{m}))^P$ in
$ (\hat{\calo}_{\hat{L}_{\epsilon}}(\widetilde{P \times_U X}^{m}))^{\glr}$
restrict to surjective morphisms
$$(\calx_{m}/\!/_{\epsilon} \CC^*)^{\hat{ss}} \to \widehat{\xu}^{m} \ \mbox{ and }
\ (\widetilde{\calx}_{m}/\!/_{\epsilon} \CC^*)^{\hat{ss}} \to \widetilde{\xu}^{m}$$
where $(\calx_{m}/\!/_{\epsilon} \CC^*)^{\hat{ss}}$ and $(\widetilde{\calx}_{m}/\!/_{\epsilon} \CC^*)^{\hat{ss}}$ 
are open subsets of
$\calx_{m}/\!/_{\epsilon} \CC^*$ and 
$\widetilde{\calx}_{m}/\!/_{\epsilon} \CC^*$
defined as in Definition \ref{defn4p4}. 
\end{theorem}

Thus when $m$ is sufficiently divisible we have the following diagrams (cf. Remark \ref{rem4p6}):
\begin{equation} \label{maindiag4} \begin{array}{ccccccccc}
(\calx_{m}/\!/_{\epsilon} \CC^*)^{\hat{ss}} &
\subseteq & \calx_{m}/\!/_{\epsilon} \CC^* & \leftarrow - \rightarrow &
X = \calx_{m}/\!/_{N/2} \CC^*  \\
 \downarrow & & & \mbox{flips} & \\
\widehat{\xu}^{m} & & & &  
\end{array}
\end{equation}
and
\begin{equation} \label{maindiag} \begin{array}{ccccccccc}
(\widetilde{\calx}_{m}/\!/_{\epsilon} \CC^*)^{\hat{ss}} &
\subseteq & \widetilde{\calx}_{m}/\!/_{\epsilon} \CC^* & \leftarrow - \rightarrow &
X = \widetilde{\calx}_{m}/\!/_{N/2} \CC^*  \\
 \downarrow & & & \mbox{flips} & \\
\widetilde{\xu}^{m} & & & &  
\end{array}
\end{equation}
where the vertical maps are {surjective} and
the inclusions are open. 
When the ring of invariants $\hat{\calo}_L(X)^U$ is finitely generated these can
be extended by Remark \ref{bbb} to 
\begin{equation} \label{maindiag3} \begin{array}{ccccccccc}
(\calx_{m}/\!/_{\epsilon} \CC^*)^{\hat{s}} & \subseteq & 
(\calx_{m}/\!/_{\epsilon} \CC^*)^{\hat{ss}} &
\subseteq & \calx_{m}/\!/_{\epsilon} \CC^* & \leftarrow - \rightarrow &
X = \calx_{m}/\!/_{N/2} \CC^*  \\
& &  
 \downarrow & & & \mbox{flips} & \\
\downarrow  & & 
\widehat{\xu}^{m} & & & &  \\
 & & \downarrow & & & & \\
X^s/U & \subseteq & 
\xu & & & &
\end{array}
\end{equation}
(and a similar diagram involving $\widetilde{\xu}^m$),
where $(\calx_{m}/\!/_{\epsilon} \CC^*)^{\hat{s}}$ is the inverse image
of $X^s/U$ in $(\calx_{m}/\!/_{\epsilon} \CC^*)^{\hat{ss}}$.
In general for sufficiently divisible $m$ and $k$ 
we have
\begin{equation} \label{maindiag5} \begin{array}{ccccccccc}
(\calx_{m}/\!/_{\epsilon} \CC^*)^{\hat{s}} & \subseteq & 
(\calx_{m}/\!/_{\epsilon} \CC^*)^{\hat{ss}} &
\subseteq & \calx_{m}/\!/_{\epsilon} \CC^* & \leftarrow - \rightarrow &
X = \calx_{m}/\!/_{N/2} \CC^*  \\
& &  
 \downarrow & & & \mbox{flips} & \\
\downarrow  & & 
\widehat{\xu}^{m} & & & &  \\
 & & \downarrow & & & & \\
X^s/U &\subseteq \xu  \subseteq & \overline{\xu}^{m'} & & & & 
\end{array}
\end{equation}
where $m'=km$.

\subsection{Naturality properties}

Given a linear action of $U = \cplusr$ on a projective variety $X$, we have  embedded
 $X$ in a normal projective variety $Y_{m}$ such that
the linear action of $U$ on $X$ extends to a linear action of $G=\slrplus \geq U$
on $Y_{m}$ and (if $m$ is sufficiently divisible) every $U$-invariant in a finite fully separating
set of $U$-invariants on $X$ extends to a $U$-invariant on $Y_{m}$. 
We have then constructed 
 the GIT quotients
$$ \widehat{\xu}^m = G \times_P (\overline{P \times_U X}^m) /\!/_{\tilde{L}_\epsilon^{(N)}} G,$$
$$\widetilde{\xu}^m = G \times_P (\widetilde{P \times_U X}^m) /\!/_{\tilde{L}_\epsilon^{(N)}} G$$
and
$$ {\calx}_m = \overline{P \times_U X}^m /\!/_{\hat{L}^{(N)}} SL(r;\CC)
\ \ \mbox{ and } \ \ \tilde{\calx}_m = \widetilde{P \times_U X}^m /\!/_{\tilde{L}^{(N)}} SL(r;\CC)$$
where $\overline{P \times_U X}^m$ and $\widetilde{P \times_U X}^m$ are the 
closures of $P\times_U X$ in the projective completions
${\prr} \times Y_{m}$ and 
$ \widetilde{\prr} \times Y_{m}$
of $P \times_U Y_{m} \cong P/U \times Y_{m}$. One difficulty with
using this construction in practice is that it is not easy to tell how divisible
$m$ has to be for there to be a finite fully separating set of $U$-invariants
on $X$ extending to $U$-invariants on $Y_{m}$. However the construction
does have the following nice property, which enables us to study
the families $\widehat{\xu}^m$ and $\widetilde{\xu}^{m}$ by embedding $X$ in other
projective varieties $Y$. 

\begin{prop} Let $Y$ be a projective variety with a very ample line bundle $L$ and 
a linear action of $U=\cplusr$ on $L$, and let $X$ be a $U$-invariant projective subvariety
of $Y$ with the inherited linear action of $U$. Then the inclusion of $X$ in $Y$
induces inclusions of projective varieties
$$\widehat{\xu}^{m} \subseteq \widehat{Y/\!/U}^{m} \ \ \mbox{ and } \ \  
\widetilde{\xu}^{m} \subseteq \widetilde{Y/\!/U}^{m}  $$
for all $m>0$, as well as ${\calx}_{m} \subseteq {\caly}_{m}$
and $\tilde{\calx}_{m} \subseteq \tilde{\caly}_{m}$ when $\calx_m$ and $\tilde{\calx}_{m}$
are defined as in Definition \ref{calxm} and $\caly_m$ and $\tilde{\caly}_{m}$ are defined similarly with
$Y$ replacing $X$.  \label{propincl}
\end{prop}

{\bf Proof}: The construction of $\widehat{\xu}^m$ and $\widetilde{\xu}^{m}$ starts by embedding $X$
into $\PP(V^X_m)$ where $V^X_m = H^0(X,L^{\otimes m})^*$,
and then embedding $\PP(V_m^X)$ into $\PP(W^X_m)$ where
$$W^X_m = \bigoplus_{j=0}^{\dim V^X_m - 1} \Theta_{j,m}^X \otimes (U_{V^X_m})^j (V^X_m)$$
and $U_{V^X_m}$ represents the infinitesimal action of $\lieu = {\rm Lie \,}(U) = \CC^r$ on
$V^X_m$. Here $\Theta^X_{j,m}$ is the space of complex valued polynomial
functions on 
$$\{ (u_1,\ldots,u_j,h_1, \ldots , h_{\dim V^X_m - 1}) \in (\CC^r)^j \times
((\CC^r)^* \otimes \CC^r)^{\dim V^X_m -1} \}$$
which are homogeneous of degree 1 in the coordinates of each $u \in \CC^r$ separately,
and homogeneous of total degree $r(\dim V^X_m - 1) - (r-1)j$ in the
coordinates of all the $h_k \in (\CC^r)^* \otimes \CC^r$ together.
The surjection $H^0(Y,L^{\otimes m}) \to H^0(X,L^{\otimes m})$ given by restriction
gives us an inclusion of $V^X_m$ into $V^Y_m$ and of $(U_{V^X_m})^j (V^X_m)$ into
$(U_{V^Y_m})^j (V^Y_m)$ for all $j \geq 0$. We then get a map
$$ \Theta^X_{j,m} \otimes (U_{V^X_m})^j (V^X_m)
\to \Theta^Y_{j,m} \otimes (U_{V^Y_m})^j (V^Y_m)$$
$$ \alpha \mapsto \alpha^Y$$
where
$\alpha^Y(u_1, \ldots,u_j,h_1, \ldots h_{\dim V^Y_m -1})$ is given by
$$ {\rm det}(h_1) \ldots {\rm det}(h_{\dim V^Y_m - \dim V^X_m})$$ times
$$\alpha(u_1,\ldots,u_j,h_{\dim V^Y_m - \dim V^X_m +1},\ldots,h_{\dim V^Y_m -1}).$$
This gives us a commutative diagram of $U$-equivariant embeddings
$$ \begin{array}{ccccccc}
X & \rightarrow & \PP(V^X_m) & \rightarrow & \PP(W^X_m) \\
\downarrow & & \downarrow & & \downarrow\\
Y & \rightarrow & \PP(V^Y_m) & \rightarrow & \PP(W^Y_m) \end{array} $$
where the righthand vertical map sends
$$\sum_{j=0}^{\dim V^X_m - 1} \alpha_j \mapsto \sum_{j=0}^{\dim V^X_m -1} \alpha_j^Y$$
and is equivariant with respect to the linear $P$-actions on $W_k^X$ and $W_k^Y$.
We thus get a commutative diagram of embeddings
$$ \begin{array}{cccccccccc}
X & \rightarrow & \PP(V^X_m) & \rightarrow &  \PP(W^X_m)
& \rightarrow & G \times_P \PP(W^X_m) \\
\downarrow & & \downarrow & & \downarrow & & \downarrow\\
Y & \rightarrow & \PP(V^Y_m) & \rightarrow &  \PP(W^Y_m) 
& \rightarrow & G \times_P \PP(W^Y_m).
\end{array} $$
$\widetilde{\xu}^{m}$ is the GIT quotient $\widetilde{G \times_U X}^{m}/\!/_{\tilde{L}_\epsilon^{(N)}}G$
where $\widetilde{G \times_U X}^{m}$ is the closure of $G(\{ \iota\} \times X) \cong G \times_U X$
in $G \times_P (\widetilde{\prr} \times (G \times_P \PP(W^X_m))$, and $\widetilde{Y/\!/U}^{m}$
is constructed in the same way. Since the line bundle $\tilde{L}_\epsilon^{(N)}$ is ample for $\epsilon >0$,
this gives us an inclusion of $\widetilde{\xu}^{m}$ into $\widetilde{Y/\!/U}^{m}$, and the other
inclusions follow similarly.

\bigskip

Suppose now that $U \cong \cplusr$ is a normal subgroup of an algebraic group
$H$ acting linearly on $X$ with respect to the line bundle $L$. This linear action induces an action of $H$
on $V_m = H^0(X,L^{\otimes m})^*$ for each $m>0$, and thus $H$ acts by conjugation on $GL(V_m)$ and its
Lie algebra. $H$ also acts by conjugation on $U$ and $\CC^r = {\rm Lie}\ U$ as $U$ is a normal subgroup of
$H$, and the Lie algebra homomorphism $\phi_{V_m} : {\rm Lie}\ U \to {\rm Lie}(GL(V_m))$
defined at (\ref{phivm}) is $H$-equivariant with respect to these actions. Hence the action of $H$ on $V_m$ 
preserves the subspaces $(U_{V_m})^j(V_m)$ of $V_m$, and from this action and the action of $H$ on $\CC^r ={\rm Lie}\ U$
we get induced actions of $H$ on $W_m$. If $h \in H$ and $p=gu \in P = \glr \ltimes U$ with $g \in \glr$ and 
$u \in U$, then the actions of $H$ and $P$ are related by
$$h(pw) = h(guw) = ((\psi(h)g\psi(h)^{-1})(huh^{-1}))(hw) = \Psi_h(p)(hw)$$
for any $w$ in $W_m$, where $\psi:H \to \glr$ is the group homomorphism defining
the action of $H$ on ${\rm Lie}\ U=\CC^r$ by conjugation and
$$\Psi_h(p) = (\psi(h)g\psi(h)^{-1})(huh^{-1}) \in \glr \ltimes U = P.$$
Thus we get actions of a semidirect product $P \rtimes H$ on $W_m$. (Note
that the action of $U$ as a subgroup of $P$ defined at (\ref{uaction})
is different from the action of $U$ as a subgroup of $H$ defined above.) In fact
the subgroup $P \rtimes U$ of $P \rtimes H$ is a direct product $P \times U$, 
since $U \cong \cplusr$ acts trivially on itself by conjugation, so if $p \in P$
and $h \in U \leq H$ then $\psi(h)$ is the identity element of $\glr$ and $\Psi_h(p)=p$.

$H$ also acts on $\prr$ and on $\widetilde{\prr}$ via the homomorphism $\psi:H \to \glr$, giving us an action of 
$P \rtimes H$ on $\widetilde{\prr} \times \PP(W_m)$. Since $X$ is $H$-invariant it follows that
the closure  $\widetilde{P \times_U X}^m$ 
  of $P \times_U X = P(\{\iota\} \times X)$ in 
$\widetilde{\prr} \times Y_m$ is also $H$-invariant. Since
$H$ normalises $SL(r;\CC)$ and commutes with the central $\CC^*$ subgroup of $\glr$,
we get an induced linear action of $H \times \CC^*$ on
$$\tilde{\calx}_m = \widetilde{P \times_U X}^m /\!/_{L^{(N)}_{N/2}} SL(r;\CC),$$
preserving the open subset
 $(\tilde{\calx}_{m}/\!/_{\epsilon} \CC^*)^{\hat{ss}}$, and of $H/U$ on 
$$\widetilde{\xu}^m = \widetilde{P \times_U X}^m /\!/_{L^{(N)}_\epsilon}P.$$

\begin{rem} \label{remcrucial}
Suppose that $\cplusr =U \unlhd H$ and $H$ acts linearly on $X$ as above.
Suppose also that $H$ contains a one-parameter subgroup $\lambda:\CC^* \to H$ whose weights for the induced
(conjugation) action on $\lieu = \CC^r$ are all strictly positive. Then 
the subgroup $\hat{U}$ of $H$ generated by $\lambda(\CC^*)$ and $U$ is a semidirect product
$$\hat{U} \cong U \rtimes \CC^*.$$
Moreover this $\CC^*$ acts on 
$\CC^{r(r+1)} \subseteq \prrplus$ with all weights strictly positive, and we have
$\glr = \psi(\lambda(\CC^*)) SL(r;\CC)$ with $\psi(\lambda(\CC^*)) \cap SL(r;\CC)$ finite.
Recall from Remark \ref{remwonderful} that if $\overline{P/U} = \widetilde{\prr}$ then 
$$\overline{P/U}^{ss,SL(r;\CC)} = \overline{P/U}^{s,SL(r;\CC)} \ \  \mbox{ and } \ \ \overline{P/U}/\!/SL(r;\CC) \cong \PP^1$$
with the induced action of $\lambda(\CC^*)$ on $\PP^1$ a positive power of the standard action on
$\CC^*$ on $\PP^1$.
Thus using variation of GIT (as in Remark 2.7) there are rational numbers
$\delta_- < \delta_+$ such that the induced action of 
$\lambda(\CC^*)$ on $\PP^1 = \overline{P/U}/\!/SL(r;\CC)$ twisted by
$\delta$ times the standard character of $\CC^*$ satisfies
$$(\PP^1)^{ss,\delta} = (\PP^1)^{s,\delta} = \left\{ \begin{array}{ccc}
\CC^* & \mbox{ if }\delta \in (\delta_-,\delta_+)\\
\emptyset & \mbox{ if }\delta \not\in [\delta_-,\delta_+] \end{array} \right. $$
and hence
$$\overline{P/U}^{ss,\glr,\delta} = \overline{P/U}^{s,\glr,\delta} = 
\left\{ \begin{array}{ccc}
\glr = P/U & \mbox{ if }\delta \in (\delta_-,\delta_+)\\
\emptyset & \mbox{ if }\delta \not\in [\delta_-,\delta_+] \end{array} \right. .
$$
It
follows that if the linearisation of the $H$-action is twisted by $\delta$ times the standard character
of $\lambda(\CC^*)$ for $\delta \neq \delta_-,\delta_+$ (which is possible to 
arrange if, for example, $\lambda(\CC^*)$ centralises $H/U$),
then (for sufficiently large $N$) all the points of $\tilde{\calx}_{m}
 = \widetilde{P \times_U X}^{m} /\!/_{L^{(N)}} SL(r;\CC)$ which are semistable for the induced
 action of $\lambda(\CC^*)$ are contained in the image of
 $ {P \times_U X} \cong \glr \times X$. Hence the inverse image
 of $(\widetilde{\xu}^{m})^{ss,\lambda(\CC^*)}$ under the surjection
 $$(\tilde{\calx}_{m}/\!/_{\epsilon} \CC^*)^{\hat{ss}} \to \widetilde{\xu}^{m}$$
 in Theorem \ref{genaction} is an open subset of $\tilde{\calx}_{m}/\!/_{\epsilon} \CC^*$
 which is unaffected by the flips $\tilde{\calx}_{m}/\!/_{\epsilon} \CC^* \leftarrow - \rightarrow
\tilde{\calx}_{m}/\!/_{N/2} \CC^* = X$ and thus can be identified canonically with an open subset
$X^{ss,\hat{U}}$ of $X$. Similarly the inverse image
 of $(\widetilde{\xu}^{m})^{s,\lambda(\CC^*)}$ under the restriction of this surjection to
 $(\calx_{m}/\!/_{\epsilon} \CC^*)^{\hat{s}}$
  is an open subset $X^{s,\hat{U}}$ of $X^{ss,\hat{U}}$ which, like $X^{ss,\hat{U}}$, is
  independent of $m$ when $m$ is sufficiently divisible. Indeed it turns out that
$x \in X^{ss,\hat{U}}$ (respectively $x \in X^{s,\hat{U}}$)
  if and only if 
$x$ is semistable (respectively stable) for every conjugate of $\lambda:\CC^* \to \hat{U}$
in $\hat{U}$, or equivalently for every one-parameter subgroup $\hat{\lambda}:\CC^* \to \hat{U}$
of $\hat{U}$. If the linear $U$-action on $X$ extends to $G$ then the same is true
when $\widetilde{\xu}^m$ is replaced with $\widetilde{\xu}$ defined as in $\S$4.1.

In particular this means that when, in addition, $H/U$ is reductive and is centralised
by $\lambda(\CC^*)$, then 
if the linearisation of the $H$-action is twisted by a suitable character of $H/U$ the induced GIT
quotients
$$  \widetilde{\xu}^{m}/\!/(H/U)
$$
for $m$ sufficiently divisible are independent of $m$
(and are isomorphic to $  \widetilde{\xu}/\!/(H/U)
$ when the linear $U$-action on $X$ extends to $G$). In fact the proof of \cite{DK} Theorem 5.3.18 shows that
in this situation, with the linearisation of the $H$-action suitably twisted, the
ring of invariants $\hat{\calo}_L(X)^H$ is finitely generated, with associated projective
variety
$$X/\!/H = {\rm Proj}( \hat{\calo}_L(X)^H),$$
and we have
$$ X/\!/H \cong \widetilde{\xu}^{m} /\!/(H/U). $$
Indeed, it turns out that in this situation, with the linearisation of the $H$-action suitably twisted,
we have essentially the same situation as for classical GIT for reductive group actions:
there is a diagram
\begin{equation} \label{maindiag2} \begin{array}{ccccccccc}
X^{{s,H}}/H & \subseteq &
X^{ss,H} & \subseteq & X  \\
\downarrow & & \downarrow & & \\
X^{s,H}/H & \subseteq & X/\!/H  & & 
\end{array}
\end{equation}
where the vertical maps are {surjective} and
the inclusions are open, and in addition the Hilbert-Mumford criteria for
stability and semistability hold as in the reductive
case (Proposition 2.3 above), and two semistable orbits
in $X$ represent the same point of $X/\!/H$ if and only if
their closures meet in $X^{ss,H}$. We will see an example of this phenomenon
for hypersurfaces in $\PP(1,1,2)$ in $\S$5.3 below.
Moreover the partial desingularisation (defined
as in $\S$2.2) $$\widetilde{\widetilde{\xu}^m/\!/(H/U)}$$ 
of the GIT quotient of $\widetilde{\xu}^{m}$ by the action of the
reductive group $H/U$ is independent of $m$ and provides
a partial desingularisation $\widetilde{X/\!/H}$ of $X/\!/H$.
\end{rem}

\begin{rem} \label{remtwisting}
In practice it is not difficult to find the range of characters 
of $\lambda(\CC^*)$ with which the linearisation can be twisted in order
to achieve the nice situation described in Remark \ref{remcrucial}. 
The picture described in Remark \ref{remcrucial} is valid for all $\delta \in \QQ \setminus \{ \delta_-,\delta_+ \}$, and if
$\delta \not\in [\delta_-,\delta_+]$ then 
$X^{ss,\hat{U},\delta} = \emptyset$ and hence $X /\!/_\delta H = \emptyset$.
If we attempt to use the Hilbert-Mumford criteria to 
calculate $X^{ss,\hat{U},\delta}$ 
 for all $\delta \in \QQ$, we will
find finitely many rational numbers $a_0< a_1 < \ldots < a_q$ such that,
when calculated according to the Hilbert-Mumford criteria,
 $X^{ss,\hat{U},\delta}$ 
  is empty for $\delta < a_0$ and for
 $\delta > a_q$, and is nonempty but constant for $\delta \in (a_{j-1},a_j)$
 when $j=1, \ldots, q$. Then we must have $\delta_- \leq a_0$ and
 $\delta_+ \geq a_q$, and moreover $X^{ss,\hat{U},\delta}$ and $X^{ss,H,\delta}$
 are as predicted by the Hilbert-Mumford criteria for any
 $\delta \neq a_0, a_q$. Thus $X/\!/_\delta H = \emptyset$ if
 $\delta \not\in [a_0,a_q]$ and the situation described in
 Remark \ref{remcrucial} holds for every $\delta \in (a_0,a_q)$.
\end{rem}

\section{Hypersurfaces in $\PP(1,1,2)$}

Recall from $\S$3 that the moduli problem of hypersurfaces
of weighted degree $d$ in the weighted projective plane $\PP(1,1,2)$ is
essentially equivalent to constructing a quotient for the action
of
$$H = (\CC^+)^3 \rtimes GL(2;\CC)$$
on (an open subset of) the projective space $X_d$ of weighted degree $d$
polynomials in the three weighted homogeneous coordinates $x,y,z$ on $\PP(1,1,2)$.
Here $H$ is the automorphism group of $\PP(1,1,2)$, where $(\a,\b,\gamma) \in
U=(\CC^+)^3$ acts on $\PP(1,1,2)$ via
$$[x:y:z] \mapsto [x:y:z + \a x^2 + \b xy + \gamma y^2]$$
and $g \in GL(2;\CC)$ acts in the standard fashion on $(x,y) \in \CC^2$
and as scalar multiplication by $(\det g)^{-1}$ on $z$. Thus $g \in GL(2;\CC)$ acts by conjugation on
$U$ as the standard action of $GL(2;\CC)$ on ${\rm Sym}^2(\CC^2) \cong \CC^3$
twisted by the character $\det$.

\begin{rem}
Notice that the central one-parameter subgroup $\lambda: \CC^* \to GL(2;\CC)$
of $GL(2;\CC)$ satisfies the conditions of Remark \ref{remcrucial} above:
the weights of its action (by conjugation) on $\lieu = \CC^3$ are all
strictly positive, as they are all equal to 4.
\end{rem}

We wish to study the action of $H$ on the projective space
$$X_d = \PP( \CC_{(d)}[x,y,z])$$
where $\CC_{(d)}[x,y,z]$ is the linear subspace of the polynomial
ring $\CC[x,y,z]$ consisting of polynomials $p$ of the form
$$p(x,y,z) = \sum_{\begin{array}{c} i,j,k \geq 0\\ i+j+2k = d \end{array}} a_{ijk}x^i y^j z^k$$
for some $a_{ijk} \in \CC$ \cite{cox,coxkatz}. Here $h \in H$ acts as $p \mapsto h\cdot p$ with
$h \cdot p(x,y,z) = p( h^{-1}x, h^{-1}y, h^{-1}z)$.
This representation $H \to GL(\CC_d[x,y,z])$ gives us a 
linearisation of the action of $H$ on $X_d$, which we can twist by any multiple
$\e \in \QQ$ of the character $\det$ of $GL(2;\CC)$ to get a fractional
linearisation $\call_\e$.

\begin{rem} \label{timesm}
If $m \geq 1$ then $H^0(X_d,\calo_{X_d}(m))^* \cong \CC_{(md)}[x,y,z]$
and the natural embedding of $X_d$ in the projective space
$\PP(H^0(X_d,\calo_{X_d}(m))^*)$ is
given by $p(x,y,z) \mapsto (p(x,y,z))^m$.
\end{rem}

\subsection{The action of $U=(\CC^+)^3$} 

First let us consider the action of the unipotent radical $U = (\CC^+)^3$
of $H$ on $X_d$. Consider
$$Y_d = \PP(\CC_{\lceil d/2 \rceil}[X,Y,W,z])$$
where $\lceil d/2 \rceil$ denotes the least integer $n \geq d/2$, and
$\CC_{\lceil d/2 \rceil}[X,Y,W,z]$ is the space of homogeneous polynomials
of degree $\lceil d/2 \rceil$ in $X,Y,W,z$. By multiplying by $x$ if $d$ is
odd, we can identify $\CC_{(d)}[x,y,z]$ with the set of polynomials of the
form
$$p(x,y,z) = \sum_{\begin{array}{c} i \geq 2\lceil d/2 \rceil - d,\ \  j,k\geq 0\\
i+j+2k = 2 \lceil d/2 \rceil \end{array} } a_{ijk}x^iy^jz^k .$$
Then we can embed $X_d$ in $Y_d$ via $p \mapsto \hat{p}$ where
$\hat{p}(X,Y,W,z)$ equals the sum over
$ i \geq 2\lceil d/2 \rceil - d$ and $ j,k\geq 0$ satisfying
$i+j+2k = 2 \lceil d/2 \rceil$ of
$$ a_{ijk} X^{(i-M_{ij})/2 - \lceil (m_{ij}-M_{ij})/2\rceil}
W^{M_{ij} + 2\lceil(m_{ij}-M_{ij})/2\rceil}Y^{(j-M_{ij})/2 - \lceil (m_{ij}-M_{ij})/2\rceil}z^k$$
for $m_{ij} = {\rm min}\{ i,j \}$ and $M_{ij} = {\rm max}\{ i,j \}$.
Thus $\hat{p}(x^2,y^2,xy,z) = p(x,y,z)$ if $d$ is even and $\hat{p}(x^2,y^2,xy,z) = xp(x,y,z)$ if $d$ is odd. For simplicity we will assume from now on that $d$ is even; by Remark \ref{timesm} this involves very little loss of generality.

The action of $U$ on $X_d$ extends to an action on $Y_d$ such that $(\a,\b,\g) \in U$ acts via
$$p(X,Y,W,z) \mapsto p(X,Y,W,z + \a X + \b W + \g Y).$$
This extends to the standard action of $G = SL(4;\CC)$ on $\CC_{d/2}[X,Y,W,z]$.
Thus 
$$Y_d/\!/U = (\PP^{12} \times Y_d)/\!/G \ \ \mbox{ and } \ \ 
\widetilde{Y_d/\!/U} = (G \times_P (\PP^{9} \times Y_d))/\!/G$$
where $\PP^{12} = \PP(\CC \oplus ((\CC^3)^* \otimes \CC^4))$ 
and $\PP^{9} = \PP(\CC \oplus ((\CC^3)^* \otimes \CC^3))$.
Here the linearisation on
$\PP^{12} \times Y_d$ is $\calo_{\PP^{12}}(N) \otimes \calo_{Y_d}(1)$
for $N>\!>0$ and 
the linearisation on
$\PP^{9} \times Y_d$ is $\calo_{\PP^{9}}(N) \otimes \calo_{Y_d}(1)$. 

The weights of the action of the standard maximal torus $T_c$ of $G=SL(4;\CC)$ on
 $\PP^{12} = \PP(\CC \oplus ((\CC^3)^* \otimes \CC^4))$ with respect to
$\calo_{\PP^{12}}(1)$ are 0 (with multiplicity 1) and $\c_1,\c_2,\c_3,\c_4$
(each with multiplicity 3) where $\c_1,\c_2,\c_3,\c_4 = -\c_1 - \c_2 - \c_3$ are the
weights of the standard representation of $SL(4;\CC)$ on $\CC^4$. The weights of the
action of $T_c$ on $Y_d = \PP(\CC_{d/2}[X,Y,W,z])$ with respect to $\calo_{Y_d}(1)$
are
$$ \{0\} \cup \{ i \c_1 + j \c_2 + k \c_3 + \ell \c_4: i,j,k,\ell \geq 0 \mbox{ and } i+j+k+\ell
= d/2 \}.$$
A point $a = [a_0:a_{11}:a_{12}:a_{13}:a_{14}:a_{21}:a_{22}:a_{23}:a_{24}:a_{31}:
a_{32}:a_{33}:a_{34}] \in \PP^{12}$ is semistable for this action of $SL(4;\CC)$
if and only if $a_0 \neq 0$. Therefore if $N>\!>0$ we have $a_0 \neq 0$ whenever
$(a,y) \in (\PP^{12} \times Y_d)^{ss,G}$ for any $y \in Y_d$. Moreover if
$a=[1:a_{ij}] \in \PP^{12}$ and $y \in Y_d$ is represented by 
$$p(X,Y,W,z) = \sum_{\begin{array}{c} i,j,k,\ell \geq 0\\i+j+k+\ell = d/2 \end{array}}
b_{ijk\ell}X^i Y^j W^k z^\ell \ \ \ \ \in \ \ \ \CC_{d/2}[X,Y,W,z]$$
then by the Hilbert-Mumford criteria $(a,y) \in (\PP^{12} \times Y_d)^{ss,G}$
if and only if $(ga,gy) \in (\PP^{12} \times Y_d)^{ss,T_c}$ for every $g \in G$,
and $(a,y) \in (\PP^{12} \times Y_d)^{ss,T_c}$ if and only if 0 lies in the convex
hull of the set of weights
$$\{ i \c_1 + j \c_2 + k \c_3 + \ell \c_4: b_{ijk\ell} \neq 0 \} \cup S_1 \cup S_2 \cup S_3 \cup
S_4$$
where
$$S_1 = \left\{  \begin{array}{cc} \{N\c_1 + i \c_1 + j \c_2 +k\c_3 + \ell \c_4:b_{ijk\ell} \neq 0
\} & \mbox{ if } (a_{11},a_{21},a_{31}) \neq 0\\
\emptyset & \mbox{ if } (a_{11},a_{21},a_{31}) = 0 \end{array} 
\right. $$
and $S_2,S_3,S_4$ are defined similarly. Let us write
\begin{equation} \label{breakup}
(\PP^{12} \times Y_d)^{ss,G} = (\PP^{12} \times Y_d)^{ss,G}_0 \sqcup
(\PP^{12} \times Y_d)^{ss,G}_1 \sqcup
(\PP^{12} \times Y_d)^{ss,G}_2 \sqcup
(\PP^{12} \times Y_d)^{ss,G}_3
\end{equation}
where $(\PP^{12} \times Y_d)^{ss,G}_q = \{(a,y) \in (\PP^{12} \times Y_d)^{ss,G}: {\rm rank}((a_{ij})) = q \}$.
Then
$$(\PP^{12} \times Y_d)^{ss,G}_0 = \{ [1:0:\ldots:0]\} \times Y^{ss,G}_d$$
and 
$$(\PP^{12} \times Y_d)^{ss,G}_1 = G \times_{U_1} (\{ [1:\iota_1]\} \times Y^{ss,1}_d)$$
where
$$ \iota_1 = \left( \begin{array}{cccc} 1&0&0&0\\ 0&0&0&0 \\ 0&0&0&0 \end{array} \right),$$
$U_1$ is its stabiliser $\{ (g_{ij} \in G: g_{11}=1 \mbox{ and } g_{21}=g_{31} = g_{41}=0 \}$
in $G = SL(4;\CC)$ and
$$Y_d^{ss,1} = \{ y \in Y_d: uy \in Y_d^{ss,T_c} \mbox{ for all }u \in U_1\}.$$
Similarly 
$$(\PP^{12} \times Y_d)^{ss,G}_2 = G \times_{U_2} (\{ [1:\iota_2]\} \times Y^{ss,2}_d)$$
where
$$ \iota_2 = \left( \begin{array}{cccc} 1&0&0&0\\ 0&1&0&0 \\ 0&0&0&0 \end{array} \right),$$
$U_2$ is its stabiliser $\{ (g_{ij} \in G: g_{11}=g_{22} =1 \mbox{ and } g_{12} = g_{21}=g_{31}=g_{32} = g_{41}=g_{42}= 0 \}$
in $G$ and
$$Y_d^{ss,2} = \{ y \in Y_d: uy \in Y_d^{ss,T^2_c} \mbox{ for all }u \in U_2\}$$
for $T^2_c = \{(g_{ij}) \in T_c: g_{11} = g_{22}\}$, while
$$(\PP^{12} \times Y_d)^{ss,G}_3 = G \times_U (\{ [1:\iota]\} \times Y^{ss,3}_d)$$
where
$$Y_d^{ss,3} = \{ y \in Y_d: uy \in Y_d^{ss,T^3_c} \mbox{ for all } u \in U\}$$
for $T^3_c = \{(g_{ij}) \in T_c: g_{11} = g_{22}=g_{33}\}$.

\subsection{The action of $\hat{U} = \CC^3 \rtimes \CC^*$} 

Now let us consider the action of the subgroup $\hat{U} = \CC^* \ltimes U$ of $H$ on
$X_d$, where $\CC^*$ is the centre of $GL(2;\CC)$ and acts by conjugation on
${\rm Lie}\,(U) = \CC^3$ with weights all equal to 4. The action of $t \in \CC^*$ on 
$p(x,y,z) \in \CC_d[x,y,z]$ is given by
$$tp(x,y,z) = p(tx,ty,t^{-2}z).$$
This action extends to the action on $\CC_{d/2}[X,Y,W,z]$ given by
$$tP(X,Y,W,z) = P(t^2X,t^2Y,t^2 W, t^{-2}z)$$
and thus the action of $\hat{U}$ on $X_d$ extends to a linear action on $Y_d$,
which is the restriction of the $GL(4;\CC)$-action on $Y_d$ via the embedding
of $\hat{U}$ in $GL(4;\CC)$ such that
\begin{equation} \label{glfour} t \mapsto \left( \begin{array}{cccc} t^2&0&0&0\\
0&t^2&0&0\\ 0&0&t^2&0\\ 0&0&0&t^{-2} \end{array} \right)     \mbox{ for }t \in \CC^*.\end{equation}
If we twist this action by $2\delta$ times the standard character of $\CC^*$ then we get a
fractional linearisation of the action of $\hat{U}$ on $Y_d$ which extends the
fractional linearisation $\call_\delta$ on $X_d$.
We also get an action of $\CC^* = \hat{U}/U$ on $\PP^{12} = \overline{G/U}$ via
$$ t[a_0:a_{ij}] = [a_0: t^4 a_{ij}].$$
Note that, since
$(t_1,t_2,t_3,t_4) = (t^2\tau_1,t^2\tau_2,t^2\tau_3,t^{-2}(\tau_1\tau_2\tau_3)^{-1})$
if and only if 
$$ t^4 = t_1t_2t_3t_4 \mbox{ and }\tau_1=t^{-2}t_1, \ \tau_2 = t^{-2} t_2, \ 
\tau_3=t^{-2} t_3,$$ 
$\CC^* T_c \cong (\CC^* \times T_c)/(\ZZ/4\ZZ)$ is the maximal torus
of $GL(4;\CC)$, acting on $\PP^{12}$ with weights $0$, $2\c_1 + \c_2 + \c_3 + \c_4$,
$\c_1 + 2\c_2 + \c_3 + \c_4$, $\c_1 + \c_2 + 2\c_3 + \c_4$ and $\c_1 + \c_2 + \c_3 + 2\c_4$,
and acting on $Y_d$ with fractional weights
$$\{ i\c_1 + j\c_2 + k \c_3 + \ell \c_4 + \frac{\e}{2}(\c_1 + \c_2 +\c_3 - \c_4): \ i,j,k,\ell \geq 0$$
$$\mbox{ and } i+j+k+\ell = \frac{d}{2}\}$$
where $\c_1,\c_2,\c_3,\c_4$ are now the weights of the standard representation of $GL(4;\CC)$ on $\CC^4$.
Let us  break up $(\PP^{12} \times Y_d)^{ss,GL(4;\CC),\delta}$
as at (\ref{breakup}) as
\begin{equation} \label{breakup2}
(\PP^{12} \times Y_d)^{ss,GL(4;\CC),\delta}_0 \ \sqcup \ 
(\PP^{12} \times Y_d)^{ss,GL(4;\CC),\delta}_1 $$
$$ \sqcup \ 
(\PP^{12} \times Y_d)^{ss,GL(4;\CC),\delta}_2 \ \sqcup \ 
(\PP^{12} \times Y_d)^{ss,GL(4;\CC),\delta}_3 
\end{equation}
where $(\PP^{12} \times Y_d)^{ss,GL(4;\CC),\delta}_q$ equals
$$ \{(a,y) \in (\PP^{12} \times Y_d)^{ss,GL(4;\CC),\delta}: {\rm rank}((a_{ij})) = q \}.$$ We find by considering the central $\CC^*$ in $GL(4;\CC)$ that
$$(\PP^{12} \times Y_d)^{ss,GL(4;\CC),\delta}_0 =
(\PP^{12} \times Y_d)^{ss,GL(4;\CC),\delta}_1 =
(\PP^{12} \times Y_d)^{ss,GL(4;\CC),\delta}_2 = \emptyset$$
unless $\delta = -d/2$, while
$$(\PP^{12} \times Y_d)^{ss,GL(4;\CC),\delta} = (\PP^{12} \times Y_d)^{ss,GL(4;\CC),\delta}_3 \cong GL(4;\CC) \times_{\hat{U}} 
Y^{ss,\hat{U},\delta}_d$$
where
$$Y_d^{ss,\hat{U},\delta} = \{ y \in Y_d: uy \in Y_d^{ss,\CC^*,\delta} \mbox{ for all }u \in U\}.$$
Similarly if $\delta \neq -d/2$ then 
$$(G \times_P (\widetilde{\PP^{9}} \times Y_d))^{ss,GL(4;\CC),\delta} = 
(G \times_P (\widetilde{\PP^{9}} \times Y_d))^{ss,GL(4;\CC),\delta}_3 $$
$$\cong GL(4;\CC) \times_{\hat{U}} 
Y^{ss,\hat{U},\delta}_d$$
and if $m$ is sufficiently divisible 
$$(G \times_P (\widetilde{P \times_U Y_d}^m))^{ss,GL(4;\CC),\delta}  \cong GL(4;\CC) \times_{\hat{U}} 
Y^{ss,\hat{U},\delta}_d,$$
while
$$(\PP^{12} \times Y_d)^{s,GL(4;\CC),\delta} 
\cong (G \times_P (\widetilde{\PP^{9}} \times Y_d))^{s,GL(4;\CC),\delta}$$
$$
\cong (G \times_P (\widetilde{P \times_U Y_d}^m))^{s,GL(4;\CC),\delta} 
\cong GL(4;\CC) \times_{\hat{U}} 
Y^{s,\hat{U},\delta}_d$$
where
$$Y_d^{s,\hat{U},\delta} = \{ y \in Y_d: uy \in Y_d^{s,\CC^*,\delta} \mbox{ for all }u \in U\}.$$
Thus if $\delta \neq -d/2$, for sufficiently divisible $m$ we have
$$\widetilde{Y_d/\!/U}^{m}/\!/_\delta \CC^* = Y_d^{ss,\hat{U},\delta}/\sim_{\hat{U}}$$
and
$$(\widetilde{Y_d/\!/U}^{m})^{s,\CC^*,\delta}/ \CC^* = Y_d^{s,\hat{U},\delta}/\hat{U},$$
and so
by Proposition 4.14 for sufficiently divisible $m$ we have
$$\widetilde{X_d/\!/U}^{m}/\!/_\delta \CC^* = X_d^{ss,\hat{U},\delta}/\sim_{\hat{U}}$$
and
$$(\widetilde{X_d/\!/U}^{m})^{s,\CC^*,\delta}/ \CC^* = X_d^{s,\hat{U},\delta}/\hat{U}$$
where
$$X_d^{ss,\hat{U},\delta} = \{ y \in X_d: uy \in X_d^{ss,\CC^*,\delta} \mbox{ for all }u \in U\}$$
and
$$X_d^{s,\hat{U},\delta} = \{ y \in X_d: uy \in X_d^{s,\CC^*,\delta} \mbox{ for all }u \in U\}$$
and $x \sim_{\hat{U}} y$ if and only if $\hat{U}x \cap \hat{U}y \cap X_d^{ss,\hat{U},\delta} \neq \emptyset$.
Using this, the proof of \cite{DK} Theorem 5.3.18 shows that
in fact, for the linearisation $\call_\delta$ when $\delta \neq -d/2$, the ring of invariants
$\hat{\calo}_{\call_\delta}(X_d)^{\hat{U}}$ is finitely generated, and
$$X_d/\!/_\delta \hat{U} = {\rm Proj}(\hat{\calo}_{\call_\delta}(X_d)^{\hat{U}}) = \widetilde{X_d/\!/U}^{m} /\!/_\delta \CC^*$$
for sufficiently divisible $m$. Thus for all $\delta \neq - d/2$
$$X_d/\!/_\delta \hat{U} = X_d^{ss,\hat{U},\delta}/\sim_{\hat{U}}$$
is a projective completion of $ X_d^{s,\hat{U},\delta}/\hat{U}$ (cf. Remark \ref{remcrucial}).

\subsection{The action of $H$} 

Let $T_c(GL(2;\CC))$ be the standard maximal torus of $GL(2;\CC) = H/U$.
It now follows immediately that 
 when $\delta \neq -d/2$, the ring of invariants
$\hat{\calo}_{\call_\delta}(X_d)^{H} = (\hat{\calo}_{\call_\delta}(X_d)^{\hat{U}})^{SL(2;\CC)}$ is finitely generated, and
\begin{equation} \label{lasteqn}
X_d/\!/_\delta H = {\rm Proj}(\hat{\calo}_{\call_\delta}(X_d)^{H}) =  X_d^{ss,H,\delta}/\sim_H \end{equation}
is a projective completion of $ X_d^{s,H,\delta}/H$, 
where
$$X_d^{ss,H,\delta} = \{ y \in X_d: uy \in X_d^{ss,GL(2;\CC),\delta} \mbox{ for all }u \in U\}$$
$$= \{ y \in X_d: hy \in X_d^{ss,T_c(GL(2;\CC)),\delta} \mbox{ for all }h \in H\}$$
and
$$X_d^{s,H,\delta} = \{ y \in X_d: uy \in X_d^{s,GL(2;\CC),\delta} \mbox{ for all }u \in U\}$$
$$= \{ y \in X_d: hy \in X_d^{s,T_c(GL(2;\CC)),\delta} \mbox{ for all }h \in H\}$$
and $x \sim_H y$ if and only if $Hx \cap Hy \cap X_d^{ss,H,\delta} \neq \emptyset$.

The weights of the action on $X_d$ of
$$T_c(GL(2;\CC)) = \left\{ \left( \begin{array}{cc} t_1 & 0 \\ 0 & t_2 \end{array} \right):t_1,t_2 \in \CC^* \right\}$$
with respect to the linearisation $\call_\delta$ are given by
$$(t_1,t_2) \mapsto t_1^{i-k+\delta} t_2^{j-k+\delta} \ \mbox{ for integers $i,j,k \geq 0$ such that $i+j+2k=d$.}$$
Thus $p(x,y,z) = \sum_{i+j+2k=d} a_{ijk}x^i y^j z^k \in \CC_{(d)}[x,y,z]$
represents a point of $X_d^{ss,T_c(GL(2;\CC),\delta}$ (respectively a point of $X_d^{s,T_c(GL(2;\CC),\delta}$)
if and only if 0 lies in the convex hull (respectively 0 lies in the interior of the convex hull) in $\RR^2$
of the subset
$$\{ (i-k+\delta,j-k+\delta): i,j,k \geq0, \ i+j+2k=d \mbox{ and } a_{ijk} \neq 0 \}$$
of the set of weights
$$\{ (i-k+\delta,j-k+\delta): i,j,k \geq0, \ i+j+2k=d \}$$
whose convex hull is the triangle in $\RR^2$ with vertices $(d+\delta,\delta)$, $(\delta,d+\delta)$ and $(\delta-d/2,\delta-d/2)$.
Notice that the bad case $\delta = -d/2$ occurs precisely when the origin lies on the edge of this
triangle joining the vertices $(d+\delta,\delta)$ and $(\delta,d+\delta)$, and that we
have $X_d^{ss,H,\delta} = \emptyset$ if $\delta \not\in (-d/2,d/2)$ (cf. Remark \ref{remtwisting}).  

Combining this with (\ref{lasteqn}) gives an explicit description of
$X_d/\!/_\delta H$ and $X^{s,H,\delta}_d/H$ whenever $\delta \neq -d/2$.

\begin{rem}
For small $d$ this description can be expressed in terms of singularities of hypersurfaces
(cf. the description in \cite{GIT} Chapter 4 $\S$2 of stability and semistability
for hypersurfaces in the ordinary projective plane $\PP^2$).
\end{rem}

\subsection{Symplectic descriptions}
 
Classical GIT quotients in complex algebraic
geometry are closely related to the process of 
reduction in symplectic geometry. 
Suppose that a compact, connected 
Lie group $K$ with Lie algebra ${\bf k}$ acts smoothly
on a symplectic manifold
$X$ and preserves the symplectic form $\omega$. 
A moment map for the action of $K$ on $X$ is a smooth map
$\mu :X\rightarrow {\bf k}^{\ast}$
which is equivariant with respect to the
given action of $K$ on $X$ and the coadjoint action of $K$ on $\lieks$,
and satisfies
$$d\mu(x)(\xi).a=\omega_x(\xi,a_x)$$
for all $x\in X$, $\xi\in T_xX$ and $a\in {\bf k}$,
where $x\mapsto a_x$ is
the vector field on $X$ defined by
the infinitesimal action of $a\in {\bf k}$.  The quotient
$\mu^{-1}(0)/K$ then inherits a symplectic structure
and is  the symplectic reduction at $0$, or
symplectic quotient, of $X$ by the action of $K$.

Now let 
$X$ be a nonsingular complex projective variety
embedded in complex projective space $\PP^n$, and let $G$
be a complex reductive group acting on $X$ via a complex linear
representation $\rho:G\rightarrow GL(n+1;\CC)$. 
If $K$ is a maximal compact subgroup of $G$, we
can choose coordinates on $\PP^n$ so that the action of $K$ preserves
the Fubini-Study form $\omega$ on $\PP^n$, which restricts to
a symplectic form on $X$. There is a moment map
$\mu :X\rightarrow {\bf k}^{\ast}$ defined by
\begin{equation} \mu(x).a = \frac{\overline{\hat{x}}^{t}\rho_{\ast}(a)\hat{x}}
{2\pi i|\!|\hat{x}|\!|^2} \label{mmap} \end{equation}
for all $a\in {\bf k}$, where $\mu(x).a$ denotes the natural pairing
between $\mu(x) \in \lieks$ and $a \in \liek$, while
$\hat{x}\in {\CC}^{n+1}-\{0\}$ is a representative
vector for $x\in \PP^n$ and the representation $\rho:K \to U(n+1)$ induces
$\rho_*: \liek \to {\bf u}(n+1)$ and dually $\rho^*:{\bf u}(n+1)^* \to \lieks$.
Note that we can think of $\mu$ as a map $\mu:X \to \liegs$ defined by
$$ \mu(x).a = \re \left( \frac{\overline{\hat{x}}^{t}\rho_{\ast}(a)\hat{x}}
{2\pi i|\!|\hat{x}|\!|^2} \right) $$
for $a \in \lieg = \liek \otimes_{\RR} \CC$; then $\mu$ satisfies
$\mu(x).a = 0$ for all $a \in i\liek$.

In this situation the GIT quotient $X/\!/G$ can be canonically
identified with the symplectic quotient 
$\mu^{-1}(0)/K$. More precisely \cite{K},
any $x\in X$ is semistable if and only if
the closure of its $G$-orbit meets $\mu^{-1}(0)$, 
while $x$ is stable if and only if its $G$-orbit
meets
$$\mu^{-1}(0)_{\rm reg} = \{ x \in \mu^{-1}(0)\  | \ d\mu(x):T_xX \to \lieks
\mbox{ is surjective} \},$$
and
the inclusions of $\mu^{-1}(0)$ into $X^{ss}$ and
of $\mu^{-1}(0)_{\rm reg}$ into $X^{s}$
induce 
homeomorphisms
$$\mu^{-1}(0)/K\rightarrow X/\!/G \mbox{ 
and }
\mu^{-1}(0)_{\rm reg} \to X^{s}/G.$$
Thus the moment map picks out a unique
$K$-orbit in each stable $G$-orbit, and also in each
equivalence class of strictly semistable
$G$-orbits, where $x$ and $y$ in $X^{ss}$
are equivalent if the closures of their
$G$-orbits meet in $X^{ss}$; that is,
if their images under the natural
surjection $q:X^{ss} \to X/\!/G$ agree.

\begin{rem}
It follows from the formula (\ref{mmap}) that if we change
the linearisation of the $G$-action of $X$ by multiplying
by a character $\chi:G \to \CC^*$ of $G$, then the
moment map is modified by the addition of a central
constant $c_\chi$ in $\lieks$, which we can identify with 
the restriction to $\liek$ of the derivative of $\chi$.
\end{rem}

When a non-reductive affine algebraic group $H$ with unipotent
radical $U$ acts linearly on a projective variety  $X$ there are
\lq moment-map-like' descriptions of suitable projectivised
quotients $\widetilde{\xu} = \widetilde{G \times_U X}/\!/G$ 
 and the resulting quotients $ 
(\widetilde{X/\!/U})/\!/(H/U)$,
which are analogous to the description of a reductive GIT quotient
$Y/\!/G$ as a symplectic quotient $\mu^{-1}(0)/K$,
and can be obtained from the symplectic quotient description of the 
reductive GIT quotient $\widetilde{G \times_U X}/\!/G$ (see \cite{K5,K6}
for more details). This is very closely
related to the \lq symplectic implosion' construction of Guillemin, Jeffrey
and Sjamaar \cite{GJS}.

The case of the automorphism group $H$ of $\PP(1,1,2)$ acting on $X_d = \PP(\CC_{(d)}[x,y,z])$
as above with respect to the linearisation $\call_\e$ for any $\e \neq -d/2$
is particularly simple. We have seen that then  $X_d/\!/_\e H = (X_d/\!/_\e \hat{U})/\!/SL(2;\CC)$
where $X_d/\!/_\e \hat{U}$ is the image of $X_d^{ss,\hat{U},\e}$ in
$Y_d/\!/_\e \hat{U} = (\PP^{12} \times Y_d)/\!/_\e GL(4;\CC)$. There is a moment
map 
$$\mu_{U(4)}: \PP^{12} \times Y_d \to {\rm Lie}\ U(4)^*$$
for the action of the maximal compact subgroup $U(4)$ of $GL(4;\CC)$ on $\PP^{12} \times Y_d$
associated to the linearisation $\calo_{\PP^{12}}(N) \times \calo_{Y_d}(1)$ for $N>\!>0$, given
by
\begin{equation} \label{ufour} \mu_{U(4)}(a,y) = N \mu_{U(4)}^{\PP^{12}}(a) \ + \ \mu_{U(4)}^{Y_d}(y)
\end{equation}
for $(a,y) \in \PP^{12} \times Y_d$. Here $\mu_{U(4)}^{\PP^{12}}:\PP^{12} \to {\rm Lie}\ U(4)^*$
and $\mu_{U(4)}^{Y_d}:Y_d \to {\rm Lie}\ U(4)^*$ are the moment maps given by formula (\ref{mmap})
for the actions of $U(4)$ on $\PP^{12}$ and on $Y_d = \PP(\CC_{d/2}[X,Y,W,z])$. We can identify
$Y_d/\!/_\e \hat{U} = (\PP^{12} \times Y_d)/\!/_\e GL(4;\CC)$ with $\mu_{U(4)}^{-1}(-\e)/U(4) = 
(\mu_{SU(4)}^{-1}(0) \cap \mu_{S^1}^{-1}(-\e))/S^1 SU(4)$, where $S^1$ is the maximal compact subgroup of
the subgroup $\CC^*$ of $GL(4;\CC)$ given at (\ref{glfour}). We can use the standard invariant inner product
on the Lie algebra of $U(4)$ to identify ${\rm Lie}\ U(4)$ with ${\rm Lie}\ U(4)^*$ and with
${\rm Lie}(U(3) \times U(1)) \oplus ({\rm Lie} (U(3) \times U(1)))^\perp$,
and thus with
$${\rm Lie} S^1 \ \oplus \ {\rm Lie}\ S(U(3)\! \times\! U(1)) \ \oplus\  ({\rm Lie}(U(3)\! \times\! U(1)))^\perp$$
where $({\rm Lie} (U(3)\! \times\! U(1)))^\perp$ is the orthogonal complement to 
${\rm Lie} (U(3)\! \times\! U(1))$ in ${\rm Lie} U(4)$, and $S(U(3)\! \times\! U(1)) = (U(3)\! \times\! U(1))\cap SU(4)$
so that $U(3)\! \times\! U(1) = S^1\ S(U(3)\! \times\! U(1))$. 
With respect to this decomposition we can write $\mu_{U(4)} = \mu_{S^1} \oplus \mu_{S(U(3) \times U(1))}
\oplus \mu_{\perp}$ where $\mu_\perp$ is the orthogonal projection of $\mu_{U(4)}$ onto
$({\rm Lie}\ (U(3) \times U(1)))^\perp$.
We find that if $\e \neq -d/2$ and $N$ is sufficiently large then
$$Y_d/\!/_{\e} \hat{U} \cong \mu_{U(4)}^{-1}(-\e)/U(4) \cong (\mu_{S^1}^{-1}(-\e) \cap \mu_\perp^{-1}(0))/S^1,$$
and restricting to $X_d$ we get an identification
$$X_d/\!/_\e \hat{U} \cong \mu_{\hat{U}}^{-1}(-\e)/S^1$$
where $\mu_{\hat{U}}:X_d \to {\rm Lie}(\hat{U})^*
\cong ({\rm Lie} S^1 \otimes_\RR \CC)^* \oplus ({\rm Lie}(U(3) \times U(1)))^\perp$
is a \lq moment map' for the action of $\hat{U}$ on
$X_d$ (which takes into account the K\"{a}hler structure on $X_d$, not just its
symplectic structure),
defined by 
$$ \mu_{\hat{U}}(x).a = \re \left( \frac{\overline{\hat{x}}^{t}\rho_{\ast}(a)\hat{x}}
{2\pi i|\!|\hat{x}|\!|^2} \right) $$
for $a \in {\rm Lie}(\hat{U})$. 
Moreover if $\e \neq -d/2$ then $X_d/\!/_\e H = (X_d/\!/_\e \hat{U})/\!/SL(2;\CC)$ can be
identified with
$$\frac{\mu^{-1}_{\hat{U}}(-\e) \cap \mu_{SU(2)}^{-1}(0)}{S^1\ SU(2)} = \frac{\mu_H^{-1}(-\epsilon)}{U(2)}$$
where the \lq moment map' $\mu_H:X_d \to {\rm Lie}(H)^*$ is 
defined by 
$$ \mu(x).a = \re \left( \frac{\overline{\hat{x}}^{t}\rho_{\ast}(a)\hat{x}}
{2\pi i|\!|\hat{x}|\!|^2} \right) $$
for $a \in {\rm Lie}(H)$, and  $U(2)$ is a maximal compact subgroup of $H$.


\begin{thebibliography}{99}

\bibitem[AD]{AD} A. Asok and B. Doran, {\it On unipotent quotients and some $\AA^1$-contractible
smooth schemes},  Int. Math. Research Papers {\bf 5} (2007), article ID rpm005.

\bibitem[AD2]{AD2} A. Asok and B. Doran, {\it Vector bundles on contractible smooth schemes},
Duke Mathematical Journal {\bf 143} (2008), 513--530.

\bibitem[AB]{AB} M.F. Atiyah and R. Bott, {\it The Yang-Mills equations over Riemann
surfaces}, Phil. Trans. Roy. Soc. London {\bf 308} (1982), 523--615.

\bibitem[BBDG]{BBDG} A. Beilinson, J. Bernstein and P. Deligne,
{\it Faisceaux pervers}, Analysis and topology on singular spaces I
(Luminy, 1981), Astérisque {\bf 100}, Soc. Math. France, Paris, 1982, 5--171.

\bibitem[BP]{BP} M. Brion and C. Procesi, {\it Action d'un tore dans une vari\'{e}t\'{e}
projective}, Progress in Mathematics {\bf 192} (1990), 509--539.

\bibitem[Co]{cox} D Cox, {\it The homogeneous coordinate ring of a toric variety},
J. Algebraic Geom. {\bf 4} (1995), 17--50. 

\bibitem[CK]{coxkatz} D Cox and S Katz, {\it Mirror symmetry and algebraic geometry}, Mathematical Surveys and Monographs {\bf 68}, American Mathematical Society, Providence, RI, 1999.

\bibitem[Do]{Dolg} I. Dolgachev, {\it Lectures on invariant theory}, London
Mathematical Society Lecture Note Series {\bf 296}, Cambridge University
Press, 2003.


\bibitem[DH]{DolgHu} I. Dolgachev and Y. Hu, {\it Variation of Geometric Invariant
Theory quotients}, Publ. Math. I.H.E.S. {\bf 87} (1998) (with an appendix by
N. Ressayre).
   

\bibitem[DK]{DK} B. Doran and F. Kirwan, {\it Towards non-reductive geometric
invariant theory}, Pure Appl. Math. Quarterly {\bf 3} (2007), 61--105.

\bibitem[DK2]{DK2} B. Doran and F. Kirwan, {\it Effective non-reductive geometric
invariant theory}, in preparation.

\bibitem[Fa]{F2} A. Fauntleroy, {\it Categorical quotients of certain algebraic
group actions}, Illinois Journal Math. {\bf 27} (1983), 115--124.

\bibitem[Fa2]{F1} A. Fauntleroy, {\it Geometric invariant theory for general algebraic
groups}, Compositio Mathematica {\bf 55} (1985), 63--87.

\bibitem[Fa3]{F3} A. Fauntleroy, {\it On the moduli of curves on rational ruled surfaces},
  Amer. J. Math.  {\bf 109}  (1987), 417--452. 

\bibitem[Gi]{Gies} D. Gieseker, {\it Geometric invariant theory and applications to moduli
problems}, { Invariant theory (Montecatini, 1982)} 
Lecture Notes in Math. {\bf 996}, 
Springer (1983), 45--73.

\bibitem[GM]{GM} M. Goresky and R. MacPherson, {\it On the topology of algebraic torus actions}, Algebraic groups, Utrecht
1986, 73--90, Lecture Notes in Mathematics, 1271. 

\bibitem[GP]{GP1} G.-M. Greuel and G. Pfister, {\it Geometric quotients of unipotent group actions},
Proc. London Math. Soc. (3) {\bf 67} (1993) 75--105.

\bibitem[GP2]{GP2} G.-M. Greuel and G. Pfister, {\it Geometric quotients of unipotent group actions II},
Singularities (Oberwolfach 1996), 27--36, Progress in Math. 162, Birkhauser, Basel 1998.

\bibitem[Gr]{Grosshans} F. Grosshans, {\it Algebraic homogeneous spaces and
invariant theory},
Lecture Notes in Math. 1673, Springer-Verlag, Berlin, 1997.


\bibitem[Gr2]{Grosshans2} F. Grosshans, {\it The invariants of unipotent radicals
of parabolic subgroups}, Invent. Math. {\bf 73} (1983), 1--9.

\bibitem[GJS]{GJS} V. Guillemin, L. Jeffrey and R. Sjamaar, {\it Symplectic implosion},
Transformation Groups {\bf 7} (2002), 155-184.



\bibitem[JK]{JK} L. Jeffrey and F. Kirwan, {\it Localization for nonabelian group
actions}, Topology {\bf 34} (1995),
291--327.

\bibitem[JKKW]{JKKW} L. Jeffrey, Y.-H. Kiem, F. Kirwan, and J. Woolf, {\it Cohomology
pairings on singular
quotients in geometric invariant theory}, Transform. Groups {\bf 8} (2003),
 217--259.

\bibitem[Ka]{Kausz} I. Kausz, {\it A modular compactification of the general
linear group}, Doc. Math. {\bf 5} (2000), 553--594.

\bibitem[KM]{MK} S. Keel and S. Mori, {\it Quotients by groupoids}, Annals
of Math. (2) {\bf 145} (1997), 193--213.



\bibitem[Ki]{K} F. Kirwan, {\it Cohomology of quotients in symplectic and
algebraic geometry},
Mathematical Notes {\bf 31} Princeton University Press, Princeton, NJ, 1984.

\bibitem[Ki2]{K2} F. Kirwan, {\it Partial desingularisations of quotients of
nonsingular varieties
and their Betti numbers}, Annals of Math. (2) {\bf 122} (1985), 41--85.

\bibitem[Ki3]{K3} F. Kirwan, {\it Rational intersection cohomology of quotient
varieties}, Invent.
Math. {\bf 86} (1986), 471--505.

\bibitem[Ki4]{K4} F. Kirwan, {\it Rational intersection cohomology of quotient
varieties II},
Invent. Math. {\bf 90} (1987), 153--167.

\bibitem[Ki5]{K5} F. Kirwan, {\it Symplectic implosion and non-reductive quotients},
to appear in the proceedings of the 65th birthday conference for Hans Duistermaat,
\lq Geometric Aspects of Analysis and Mechanics', Utrecht 2007.

\bibitem[Ki6]{K6} F. Kirwan, {\it Generalised symplectic implosion}, in preparation.

\bibitem[Muk]{Mukai1} S. Mukai, {\it An introduction to invariants and moduli},
Cambridge University Press 2003.

\bibitem[Muk2]{Mukai} S. Mukai, {\it Geometric realization of $T$-shaped root systems
and counterexamples
to Hilbert's fourteenth problem}, Algebraic transformation groups and
algebraic varieties, 123--129, Encyclopaedia Math. Sci. {\bf 132},
Springer, Berlin, 2004.


\bibitem[MFK]{GIT} D. Mumford, J. Fogarty and F. Kirwan, {\it Geometric invariant
theory}, 3rd edition, Springer, 1994.




\bibitem[Na]{Nagata} M. Nagata, {\it On the $14$-th problem of Hilbert},
Amer. J. Math. 81, 1959, 766--772.



\bibitem[Ne]{New} P.E. Newstead, {\it Introduction to moduli problems and orbit
spaces}, Tata Institute Lecture Notes, Springer, 1978.



\bibitem[Ne2]{New2} P.E. Newstead, {\it Geometric invariant theory}, these proceedings.



\bibitem[Po]{Popov} V. Popov, {\it On Hilbert's theorem on invariants},
Dokl. Akad. Nauk SSSR 249 (1979), 551--555.  English
translation: Soviet Math. Dokl. 20 (1979), 1318--1322 (1980).


\bibitem[PV]{PopVin} V. Popov and E. Vinberg, {\it Invariant theory},
Algebraic geometry IV, Encyclopaedia of Mathematical Sciences v. 55,
1994.

\bibitem[Rei]{Reichstein} Z. Reichstein, {\it Stability and equivariant maps}, Invent. Math. {\bf 96}
(1989), 349--383.

\bibitem[Res]{Ress} N. Ressayre, {\it The GIT-equivalence for $G$-line bundles}, Geom. Dedicata {\bf 81}
(2000), 295--324.


\bibitem[Th]{Thad} M. Thaddeus, {\it Geometric invariant theory and flips}, Journal of
Amer. Math. Soc. {\bf 9} (1996), 691--723.

\bibitem[Wi]{W} J. Winkelmann, {\it Invariant rings and quasiaffine quotients}, Math. Z. 244 (2003), 163--174.



\end{thebibliography}
\end{document}